\documentclass[12pt]{article}

\usepackage{amssymb}
\usepackage{amsmath}
\usepackage{amsbsy}
\usepackage{amscd}
\usepackage{amsfonts}
\usepackage{amsthm}
\usepackage{mathrsfs}
\usepackage{verbatim}
\usepackage{hyperref}
\usepackage{fullpage}
\usepackage{mathdots}
\usepackage{graphicx,subfigure}
\usepackage[english]{babel}
\usepackage[utf8]{inputenc}
\usepackage{tikz}
\usepackage{adjustbox}
\usepackage{authblk}
\usepackage{caption}

\usepackage{mathtext}

\usepackage{tikz-cd}
\usetikzlibrary{backgrounds,fit, matrix}
\usetikzlibrary{positioning}
\usetikzlibrary{calc,through,chains}
\usetikzlibrary{arrows,shapes,snakes,automata, petri}
\usepackage[boxsize=1.25em, centerboxes]{ytableau}

\newtheorem{Theorem}[equation]{Theorem}
\newtheorem{Corollary}[equation]{Corollary}
\newtheorem{Lemma}[equation]{Lemma}
\newtheorem{Proposition}[equation]{Proposition}

\theoremstyle{definition}
\newtheorem{Definition}[equation]{Definition}
\newtheorem{Example}[equation]{Example}

\newtheorem{Notation}[equation]{Notation}

\newtheorem{Remark}[equation]{Remark}

\numberwithin{equation}{section}
\numberwithin{figure}{section}

\newcommand{\C}{{\mathbb C}}
\newcommand{\Z}{{\mathbb Z}}
\newcommand{\Q}{{\mathbb Q}}

\newcommand{\mc}[1]{\mathcal{#1}}

\newcommand{\mb}[1]{\mathbb{#1}}

\DeclareMathOperator{\MO}{\textbf{\"O}}

\begin{document}

	\title{Applications of Homogeneous Fiber Bundles to the Schubert Varieties}
	
	\author[1]{Mahir Bilen Can}
	\author[2]{Pinaki Saha}
	
	\affil[1]{\small{Tulane University, New Orleans, Louisiana\\mahirbilencan@gmail.com}} 
	\affil[2]{\small{Indian Institute of Technology Bombay, Mumbai, pinaki653@gmail.com}}

\maketitle

\begin{abstract}
This article explores the relationship between Schubert varieties and equivariant embeddings, using the framework of homogeneous fiber bundles over flag varieties. 
We show that the homogenous fiber bundles obtained from Bott-Samelson-Demazure-Hansen varieties are always toroidal.
Furthermore, we identify the wonderful varieties among them. 
We give a short proof of a conjecture of Gao, Hodges, and Yong for deciding when a Schubert variety is spherical with respect to an action of a Levi subgroup. 
By using BP-decompositions, we obtain a characterization of the smooth spherical Schubert varieties. Among the other applications of our results are: 1) a characterization of the spherical Bott-Samelson-Demazure-Hansen varieties, 
2) an alternative proof of the fact that, in type A, every singular Schubert variety of torus complexity 1 is a spherical Schubert variety, and 3) a proof of the fact that, for simply laced algebraic groups of adjoint type, every spherical $G$-Schubert variety is locally rigid,
that is to say, the first cohomology of its tangent sheaf vanishes.



\end{abstract}

\section{Introduction}

Let $G$ be a connected reductive group defined over an algebraically closed field. 
Let $X$ be a normal $G$-variety. 
If a Borel subgroup $B$ of $G$ has an open orbit in $X$, then $X$ is said to be a {\em spherical $G$-variety}. 
This conventional definition of a spherical variety requires a reductive group action, but it is worth recognizing that there are some outstanding examples of non-$G$ but $B$-varieties $X$ that have an open orbit for $B$. Some of the most prominent instances of such varieties include Schubert varieties and some Bott-Samelson-Demazure-Hansen varieties (BSDH-varieties). Moreover, such varieties always have a reductive group action with its Borel subgroup contained in $B$.
So we let $X$ denote a Schubert variety in the full flag variety $G/B$. 
In the present article, we investigate, among other things, the following basic question: 
	\medskip
	
	Is $X$ a spherical $L$-variety, where $L$ is a maximal reductive subgroup of the stabilizer of $X$ in $G$?
	\medskip

	Let $T$ be a maximal torus contained in $B$, and let $W$ denote the Weyl group of the pair $(G,T)$. 
	Let $S$ be the set of simple reflections of $W$ relative to $B$. 
	For $I\subset S$, we denote by $W_I$ (resp. by $W^I$) the subgroup generated by $I$ in $W$ (resp. 
	the set of minimal length left coset representatives of $W_I$ in $W$). 
	Let $P_I$ denote the parabolic subgroup generated by $B$ and the representatives in $G$ of the elements of $W_I$.
	For $w\in W^I$, the associated {\em Schubert variety} in $G/P_I$, denoted by $X_{wP_I}$, 
	is defined as the Zariski closure of the $B$-orbit of the point $wP_I$ in $G/P_I$. 
	In particular, since $W^\emptyset = W$, the Schubert varieties in $G/B$ are indexed by the elements of $W$.
	Let $w\in W$. 
	The stabilizer of the Schubert variety $X_{wB}$ is a parabolic subgroup $P(w)$ such that $B\subset P(w)$.
	Let $L(w)$ be a Levi subgroup of $P(w)$ such that $T\subset L(w)$. 
	Then the Weyl group of $L(w)$ is of the form $W_{J(w)}$ for some $J(w)\subseteq S$. 
	For $J\subset J(w)$, we denote by $W_J, L_J$, and $B_J$, respectively, the Weyl group generated by $J$,
	the Levi subgroup of $L(w)$ determined by $W_J$, and the Borel subgroup $B\cap L_J$. 
	By a {\em Coxeter element} we mean the product of all elements of $S$ in some order. 
	In this terminology, the first main result of our paper is the following statement.
	
\begin{Theorem}\label{T:major}
Let $G$ be a connected semisimple simply connected algebraic group defined over an algebraically closed field. 
Let $X_{wB}$ be a Schubert variety in $G/B$.
Then, for $J\subseteq J(w)$, the following statements are equivalent: 
\begin{enumerate}
\item $X_{wB}$ is an $L_J$-spherical variety such that $\dim B_J= \dim X_{wB}$;
\item $w=w_{0,J}c$, where $w_{0,J}$ is the longest element of $W_{J}$ and $c$ is a Coxeter element of $W$ such that $\ell(w)=\ell(w_{0,J})+\ell(c)$. 
\end{enumerate}
\end{Theorem}

The statement of our theorem was proposed as a conjecture by Gao, Hodges, and Yong in~\cite{GHY}, where the underlying field is assumed to be of characteristic 0, and $G$ was assumed to be simple. 
In the original version of the Hodges-Yong conjecture, the equality $\dim B_J= \dim X_{wB}$ was not required.  
By adopting the proof of our theorem to the action of a torus quotient of $B_J$, we easily remove this condition. 
We should mention that, in type A, the Hodges-Yong conjecture has recently been proved in~\cite{gao2021classification} by Gao, Hodges, and Yong.\footnote{After the completion of this article, we learned that 	Gao, Hodges, and Yong~\cite{GaoHodgesYong2023} independently proved their conjecture.} 
Working with the Demazure characters in characteristic zero, they use a totally different, computational approach.
We should also mention that, in the same article~\cite{gao2021classification}, the authors conjectured a pattern avoidance criterion for the sphericalness of a Schubert variety in type A. More recently, in~\cite[Theorem 1.4]{Gaetz}, Gaetz proved this pattern-avoidance criterion. 
\medskip
	
As an application of Theorem~\ref{T:major} and a result of Richmond and Slofstra from~\cite{RichmondSlofstra}, we prove the following interesting theorem which shows that there are strong connections between smoothness and sphericalness of Schubert varieties. 
	
\begin{Theorem}\label{T:inverseissmooth}
Let $w\in W$. 
Let $J$ be a subset of $J(w)$.
We assume that $X_{wB}$ is a spherical $L_J$-variety such that $\dim X_{wB}=\dim B_{L_J}$.  
Then the following assertions are equivalent: 
\begin{enumerate}
\item $X_{w B}$ is a smooth Schubert variety.
\item $X_{w^{-1}B}$ is a smooth Schubert variety.
\item $X_{c^{-1}P_J}$ is a smooth toric variety.
\end{enumerate}
\end{Theorem}

Note that the equivalence of the first two assertions in the theorem was proven by Carrell~\cite{Carrell2011}.

	\medskip

Our next result establishes a connection between the BSDH-varieties and equivariant embeddings. 
Let $\underline{w}$ be a word in the simple reflections of $W$ relative to $B$. 
We denote by $X_{\underline{w}}$ the BSDH-variety associated with $\underline{w}$. 
It is defined as follows. 
		If $\underline{w} = (s_{i_1},\dots, s_{i_m})$, then we have the product of minimal parabolic subgroups $\prod_{j=1}^m P_{i_j}$,
		where $P_{i_j} = B\cup B s_{i_j} B$ for $j\in \{1,\dots, m\}$.
		The corresponding {\em BSDH-variety}, denoted $X_{\underline{w}}$, is defined as the quotient variety $\prod_{j=1}^m P_{i_j} / B^m$, 
		where the (right) action of $B^m$ is given by 
		\[
		(p_1,\dots, p_m) \cdot (b_1,\dots, b_m) = (p_1 b_1, b_1^{-1} p_2 b_2 , \dots, b_{m-1}^{-1} p_m b_m),
		\]
		for $(b_1,\dots, b_m)\in B^m$ and $(p_1,\dots, p_m)\in P_{i_1}\times \cdots \times P_{i_m}$. 
		
Let $w$ be the element of $W$ that is obtained by multiplying the entries of $\underline{w}$ in the order that they appear in $\underline{w}$. 
If $\underline{w}$ is a reduced word, then the length of $\underline{w}$ is equal to the dimension of $X_{wB}$. 
Furthermore, in this case, there is a natural map $\mathbf{m}: X_{\underline{w}} \to X_{wB}$ that is a resolution of singularities of $X_{wB}$. 
The homogenous fiber bundles of the form $G\times_B X_{wB}$ are known as the $G$-Schubert varieties (in $G/B\times G/B$). 
In this paper, we term a homogeneous fiber bundle of the form $G\times_B X_{\underline{w}}$ a {\em $G$-BSDH variety}. 
Note that, since the BSDH-varieties are nonsingular, so are the $G$-BSDH varieties. 
Indeed, a $G$-BSDH variety can be viewed as a natural $G$-equivariant resolution of singularities of a $G$-Schubert variety via the canonical map 
$\mathbf{m}_G := \mathbf{1}\times \mathbf{m} : G\times_B X_{\underline{w}} \to G\times_B X_{wB}$, $[g,x]\to [g, \mathbf{m}(x)]$.
Now let $H$ be a closed subgroup of $G$ such that the left multiplication action of $B$ on $G/H$ has an open orbit. 
Then $H$ (resp. $G/H$) is called a {\em spherical subgroup} (resp. a {\em spherical homogeneous space}). 
A $G$-equivariant embedding $X$ of $G/H$ is said to be {\em simple} if $G$ has a unique closed orbit in $X$.
It is called a {\em toroidal embedding} if whenever a $B$-stable prime divisor $D$ of $X$ contains a $G$-orbit, then $D$ is $G$-stable. 
Our second main result is the following statement.

\begin{Theorem}\label{T:toroidal}
Let $X$ be a $G$-BSDH variety.
Let $D$ be a $B$-stable divisor in $X$. 
If $D$ contains a $G$-orbit, then $D$ is $G$-stable. 
\end{Theorem}
	
This result shows that a $G$-BSDH variety behaves like a toroidal variety even if it is not spherical.
By using Theorem~\ref{T:toroidal}, we obtain a similar result for $G$-Schubert varieties. 

\begin{Corollary}\label{C:toroidalforGSchubert}
Let $G\times_B X_{wB}$ be a $G$-Schubert variety.
Let $D$ be a $B$-stable divisor in $G\times_B X_{wB}$. 
If $D$ contains a $G$-orbit, then $D$ is $G$-stable. 
\end{Corollary}

As an example of such a $G$-Schubert variety, we may consider $X:=G\times_B X_{w_0 B}$, where $w_0$ is the longest element of $W$.
Then we have a $G$-equivariant isomorphism $X\to G/B \times G/B$.
On the target, the action of $G$ is the diagonal action. 
The classification of spherical diagonal actions on double partial flag varieties is well-known. 
In~\cite{Littelmann}, Littelmann classified the spherical diagonal actions on the products of two Grassmannians. 
Working in type A, Magyar, Weyman, and Zelevinsky~\cite{MWZ1} classified the spherical diagonal actions across all products of partial flag varieties.
Finally, Stembridge classified spherical diagonal actions on double partial flag varieties for all types in \cite{Stembridge}.
In particular, we know from~\cite[Theorem 2.2]{MWZ1} (or from~\cite[Corollary 1.3.A]{Stembridge}) that, for $G=GL(n,\C)$, where $n\geq 4$, $G/B\times G/B$ is not a spherical $G$-variety. Therefore, $X$ cannot be a spherical $G$-variety. 
Nevertheless, Corollary~\ref{C:toroidalforGSchubert} implies that if a $B$-stable divisor $D$ in $X$ contains a $G$-orbit, then $D$ is a $G$-stable divisor.  
\medskip

We proceed with the assumption that $G/H$ is a spherical homogeneous space as before. 
A $G$-equivariant embedding $X$ of $G/H$ is called a {\em wonderful variety} if it is smooth, complete, simple, and toroidal.  
Strictly speaking, the original definition of a ``wonderful variety'' is slightly different. 
The equivalence of the definitions is the main result of Luna's paper~\cite{Luna1996}.
By building on some basic observations of Luna and Avdeev, we obtain the following interesting connection between Schubert varieties and wonderful varieties. 

\begin{Theorem}\label{T:main2}
Let $\underline{w}$ be a word in $S$. Then $G\times_B X_{\underline{w}}$ is a wonderful variety if and only if $X_{\underline{w}}$ is a toric variety. 
Furthermore, if $\underline{w}$ is a reduced word in $S$, then $G\times_B X_{\underline{w}}$ is a wonderful variety if and only if $X_{wB}$ is a toric variety. 
\end{Theorem}

The structure of our paper is as follows. 
In the next section, we review some basic invariants of the algebraic group actions. 
In Section~\ref{S:Gorbs}, we investigate the poset of $G$-orbit closures in a homogenous fiber bundle over $G/H$. 
We prove the following fact (Theorem~\ref{T:GorbitinterectsY}): the inclusion poset of $H$-orbit closures in the fiber over $eH$ is isomorphic to the inclusion poset of $G$-orbit closures in $X$, where $X$ is assumed to possess an open $G$-orbit. 
In Section~\ref{S:GBSDH}, we prove Theorem~\ref{T:toroidal}.
In Section~\ref{S:Toroidal}, we prove Theorem~\ref{T:main2}.
In Section~\ref{S:PD}, we prove our first main result, Theorem~\ref{T:major}.
We discuss some applications of our results in Section~\ref{S:Applications}.
In particular, we prove Theorem~\ref{T:inverseissmooth} in this final section. 
We finish our paper by mentioning some future work on the relative versions of our results.

\section{Preliminaries}\label{S:Prelims}

Throughout this article, we work over an algebraically closed field $k$. 
The one-dimensional multiplicative (resp. additive) algebraic groups $(k^\times, \cdot)$ (resp. $(k,+)$) will be denoted by $\mathbb{G}_m$ (resp. by $\mathbb{G}_a$).
The letter $G$ is reserved for a connected affine algebraic group unless otherwise specified.
We will use the letter $B$ to denote a Borel subgroup of $G$. 
Unless otherwise noted, the letter $e$ will stand for the identity element of an algebraic group under consideration. 
\medskip

	If an algebraic group $H$ acts by a morphism on an algebraic variety $X$, then we will denote the action by $H:X$.

\subsection{Modality and complexity.}

Let $X$ be a $G$-variety. 
The {\em generic modality} of the action $G:X$, denoted by $d_G(X)$, is the transcendence degree over $k$ of the field of $G$-invariant rational functions on $X$. 	
In other words, we have $d_G(X) := \text{tr.deg}\ k(X)^G$.
It follows from a well-known result of Rosenlicht~\cite[Theorem 2]{Rosenlicht1956} on the rational invariants that the generic modality is equal to the minimum codimension of a $G$-orbit in $X$. A readily accessible proof of this fact can be found in~\cite[Corollary 2.3]{PopovVinberg1994}.
The {\em modality} of $G:X$, denoted by $\text{mod}(G:X)$, is defined by $\text{mod}(G:X) = \max_{Y\subseteq X} d_G(Y)$, where $Y$ varies in the set of all $G$-stable irreducible subvarieties of $X$. 
It follows from a result of Popov and Vinberg in~\cite[Theorem 8]{PopovVinberg1972} that the actions with modality 0 are precisely the actions with a finite number of orbits.
	
	\medskip
	For a reductive group $G$, the {\em complexity} of the action $G:X$, denoted by $c_G(X)$, can be defined as the generic modality of $X$ with respect to the action $B:X$, where $B\subset G$ is a Borel subgroup. 
	In particular, the condition $c_G(X)=0$ is equivalent to $B$ having an open orbit in $X$. 
	If $X$ is normal, $G$ is reductive, and $c_G(X) = 0$, then $X$ will be called a {\em spherical $G$-variety}. 
	Brion and Vinbeg showed separately that $X$ is a $G$-spherical variety if and only if there are only finitely many $B$-orbits in $X$.
	In other words, $\text{mod}(B:X)=0$ if and only if $c_G(X)=0$.
	\medskip

	Let $H$ be a closed subgroup of $G$. 
	If $G/H$ is a spherical $G$-variety via the left multiplication action, then $H$ is called a {\em spherical subgroup}.
	In this case, we will refer to $G/H$ as a {\em spherical homogeneous space}.
	If there is an open $G$-orbit $O$ in a normal $G$-variety $X$, then we will call $X$ an {\em embedding of $O$}.
	In particular, any spherical $G$-variety is an embedding of a spherical homogeneous space.

\subsection{Toroidal embeddings, regular varieties.}\label{SS:toroidalregular}

Let $G$ be a connected reductive group. 
Let $X$ be a spherical $G$-variety.  
Let $H$ denote the stabilizer in $G$ of a point $x_0$ from the open $G$-orbit in $X$. 
Then $G/H \cong G\cdot x_0$. 
Let us denote the open orbit by $X_0$. 
We denote the set of $B$-stable prime divisors of $X_0$ by $\mathcal{D}$. 
Notice that if a prime divisor $D$ does not intersect $X_0$, then $D$ must be $G$-stable. 
Indeed, the closed set $X\setminus X_0$ is $G$-stable, and $D$ is one of its irreducible components.
Since $G$ is connected, the closures of its orbits are the irreducible components of $X\setminus X_0$. 
Notice also that a $G$-stable prime divisor $E\subseteq X$ cannot intersect the open orbit nontrivially. 
Otherwise, we would have $E=X$, which is absurd. 
Thus, we see that the $G$-stable prime divisors of $X$ are precisely the irreducible components of $X\setminus X_0$.
		
	Let $Y \subset X$ be a $G$-orbit. 
	We set 
	\begin{align}\label{A:notationtoroidal}
	\mc{D}_Y := \{ D\in \mc{D} :\ Y\subseteq \overline{D} \}\qquad\text{and}\qquad \mc{D}_{all}:=\bigcup_{\text{$Y$ is a $G$-orbit in $X$}} \mc{D}_Y.
	\end{align}
	A $B$-stable but not $G$-stable divisor of $X$ is called a {\em color of $X$}. 
	In other words, the colors of $X$ are the Zariski closures in $X$ of the elements of $\mc{D}$. 
	A spherical variety $X$ is called {\em toroidal} if none of its colors contain a $G$-orbit. 
	Equivalently, $X$ is toroidal if $\mc{D}_{all} = \emptyset$.
	In this case, we call $X$ a {\em toroidal embedding} of $G/H$. 
	In our earlier work,~\cite{CanHodgesLakshmibai}, we determined some necessary conditions for a Schubert $X_{wB}$ to be toroidal.
	\medskip

	The notion of a ``regular $G$-variety'' is introduced by Bifet, De Concini, and Procesi in~\cite[Definition 5]{BDP}. 
	A smooth $G$-variety $X$ is called a {\em regular $G$-variety} if the following three conditions are satisfied:
	\begin{enumerate}
		\item $X$ contains an open $G$-orbit $X_0$ such that $X\setminus X_0$ is a union of smooth prime divisors with normal crossings.
		These prime divisors are called the {\em boundary divisors}. 
		\item Every $G$-orbit closure in $X$ is the transversal intersection of the boundary divisors. 
		\item For every $x\in X$, the normal space $T_x X / T_x (G\cdot x)$ contains a dense orbit of ${\rm Stab}_G(x)$. 
	\end{enumerate}
	It turns out that any complete regular $G$-variety is spherical. 
	Conversely, every homogeneous spherical $G$-variety admits a completion, which is a regular $G$-variety. 
	\medskip
	
	Under the assumption of completeness, the toroidal and regular varieties are closely related to one another. 
	The following result which we will use later is due to Bien and Brion~\cite[Proposition 2.2.1]{BienBrion}.
	\begin{Theorem}\label{T:BienBrion}
	Let $X$ be a smooth complete spherical $G$-variety. Then $X$ is a toroidal embedding (of its open orbit) if and only if $X$ is a regular $G$-variety. 
	\end{Theorem}

\subsection{Homogeneous fiber bundles.}\label{SS:Homogeneous}
	
Let $H$ be a closed subgroup of $G$. 
The notion of a ``homogeneous bundle'' over $G/H$ is a bridge between the category of $H$-varieties and the category of $G$-varieties.
A concise but good presentation of this useful gadget is given in~\cite[Chapter 2.1]{Timashev}.
\begin{Definition}
A {\em homogeneous fiber bundle over $G/H$} is a $G$-variety $X$ together with a $G$-equivariant surjective morphism $X\to G/H$.
The homogeneous space $G/H$ is called the {\em base} of homogeneous fiber bundle. 
\end{Definition}

Quotients provide for a more precise expression of homogenous fiber bundles.
Let $Z$ be a quasi-projective $H$-variety. 
Then $H$ acts diagonally on $G\times Z$ via $h\cdot (g,z) := (gh^{-1}, h\cdot z)$ for $h\in H$, and $(g,z)\in G\times Z$. 
The quotient set, denoted $G\times_H Z$, is a $G$-variety. 
In fact, $G\times_H Z$ is a homogeneous fiber bundle over $G/H$; the surjective $G$-equivariant morphism is given by 
\begin{align*}
G\times_H Z &\longrightarrow G/H\\
[(g,z)] &\longmapsto gH.
\end{align*} 	
Some simple examples of homogeneous fiber bundles will be useful for our purposes.

\begin{Example}\label{E:BWB}
Let $G$ be a connected reductive group.
Let $B$ be a Borel subgroup of $G$. 
Let $\chi$ be a character of $B$.
Then we have a one dimensional representation $\rho :B\to \mb{G}_m$ defined by $\rho(b) \cdot x := \chi(b) x \ (x \in k,\ b\in B)$.
By using this action of $B$ on $k$, we get a $B$-variety $G\times k$ where the action of $B$ is given by 
\begin{align}\label{A:BorelWeilBott}
b\cdot (g,x) = (gb^{-1}, \rho(b)\cdot x) = (gb^{-1}, \chi(b)x)\qquad (g\in G, x\in k, b\in B).
\end{align}
The action (\ref{A:BorelWeilBott}) has a geometric quotient.
Hence, its quotient is an algebraic variety, which we denote by $G \times_B k_\chi$. 
Clearly, the map 
\begin{align*}
p: G\times_B k_\chi &\longrightarrow G/B \\
[g,x] &\longmapsto gB
\end{align*}
is a surjective morphism.
Furthermore, $G$ acts on $G\times_B k_\chi$ via left multiplication on the first factor. 
The projection $p: G\times_B k_\chi \to G/B$ is equivariant with respect to this action. 
In other words, $G\times_B k_\chi$ is a homogeneous fiber bundle over $G/B$. 
It is easy to check that $G\times_B k_\chi$ is (the total space of) a line bundle, denoted by $\mathcal{L}_\chi$, on the flag variety $G/B$. 
In characteristic 0, if $\chi$ is an anti-dominant weight of $B$, then the dual of the space of global sections of $\mathcal{L}_{\chi}$, that is $H^0(G/B, \mathcal{L}_\chi)^*$, is the irreducible representation of $G$ with highest weight $\chi$. 
For further details of this correspondence, see~\cite[Chapter II]{Jantzen}.
\end{Example}

	The previous example has a far reaching generalization. 
	Let $P$ be a parabolic subgroup of $G$.
	Let $Y$ be a $P$-variety. 
	Then $G\times_P Y \to G/P$ is a homogeneous fiber bundle with fiber $Y$. 
	In this case, there is an equivalence between the category of $G$-linearized sheaves on the homogeneous fiber bundle $G\times_P Y$ and the category of $P$-linearized sheaves on $Y$, see~\cite[Section 2]{Brion2003} and~\cite[Section 4]{HeThomsen2008}.

	\medskip

	We continue with a very simple lemma that has some important consequences.
	\begin{Lemma}\label{L:restrictionishomogeneous}
		Let $X$ be a homogeneous fiber bundle over $G/H$. 
		If $O$ is a $G$-orbit (closure) in $X$, then $O$ is a homogeneous fiber bundle over $G/H$.
	\end{Lemma}
	\begin{proof}
		We will show that the restriction of $p$ to $O$, denoted $p|_O$, is a $G$-equivariant surjective morphism. 
		Let $z$ be a point from $O$. 
		Then, by the $G$-equivariance of the original map $X\to G/H$, we have $p ( G\cdot z ) = G \cdot p(z) = G/H$. 
		It follows that $p |_O$ is surjective. 
		This finishes the proof of our assertion. 
	\qed\end{proof}

The structure of a $G$-orbit in a homogeneous fiber bundle is determined to a certain degree by the fiber at the ``origin'' of the base. 
	
\begin{Lemma}\label{L:GorbitinterectsY1}
Let $p: X\to G/H$ be a homogeneous fiber bundle over $G/H$. 
Let $Y$ denote the fiber at $eH$, that is, $Y:=p^{-1}(eH)$.
Then every $G$-orbit in $X$ intersects $Y$. 
Furthermore, $Y$ is stable under $H$-action. 
\end{Lemma}
\begin{proof}
Let $x$ be a point in $X$. 
Since $p$ is a $G$-equivariant morphism, we have $p(G\cdot x) = G\cdot p(x)$. 
Clearly, $G\cdot p(x)$ is equal to $G/H$.
In particular, there exists $g\in G$ such that $p(g\cdot x) = eH$.
Therefore, we have $g\cdot x \in (G\cdot x) \cap Y$.
This means that the intersection $(G\cdot x) \cap Y$ is nonempty.  

Our second claim also follows from the $G$- (hence $H$-) equivariance of $p$ combined with the fact that $HeH=H=eH$.
\qed\end{proof}

\begin{Lemma}\label{L:GorbitinterectsY2}
We maintain the notation from Lemma~\ref{L:GorbitinterectsY1}.
Let $g\in G$ and $y\in Y$. 
Then we have $g\cdot y \in Y$ if and only if $g\in H$. 
\end{Lemma}
\begin{proof}
If $g\cdot y \in Y$, then we have $p(g\cdot y) = eH$.
By the $G$-equivariance of $p$, we see that $p(g\cdot y)  = g\cdot p(y)$. 
But $y\in Y$ implies that $p(y) = eH$, hence that, $gH = eH$, or $g\in H$. 
Conversely, if $g\in H$, then $p(g\cdot y) = g\cdot eH = gH = eH$. 
Therefore, $g\cdot y \in Y$. This finishes the proof of our assertion. 
\qed\end{proof}

\begin{Corollary}\label{C:Gorbs}
We maintain the notation from Lemma~\ref{L:GorbitinterectsY1}.
We assume that $Y$, that is, the fiber at $eH$, is irreducible. 
Then $G$ has an open orbit in $X$ if and only if $H$ has an open orbit in $Y$.
\end{Corollary}	

\begin{proof}
By Lemmas~\ref{L:GorbitinterectsY1} and~\ref{L:GorbitinterectsY2}, we see that 1) $Y$ is $H$-stable, and 2) every $G$-orbit in $X$ intersects $Y$ along an $H$-orbit in $Y$. 
It follows that if $O$ is an open $G$-orbit in $X$, then $O\cap Y$ is an open $H$-orbit in $Y$.
To prove the converse statement, let us assume that $H$ has an open orbit, denoted by $O$, in $Y$. 
Towards a contradiction, let us assume also that $G$ does not have an open orbit in $X$. 
Then there are infinitely many $G$-orbits of maximal dimension in $X$. 
This follows from the well-known fact that~\cite[Ch. 7, Theorem 3.3]{FerrerSantosRittatore} the points whose orbits have maximal dimension form an open subset $C\subset X$.
Let us show that $C\cap O$ is a nonempty open subset of $Y$. 
Indeed, $C\cap Y$ is a nonempty open subset of $Y$. 
Since $Y$ is irreducible, all open subsets in $Y$ intersect each other. 
In particular, we see that $C\cap O\neq \emptyset$. 
Now, since $C\cap O$ contains infinitely many points with disjoint $G$-orbits, we find a contradiction. 
Hence, we conclude that $G$ has an open orbit in $X$.
\qed\end{proof}

\section{The Poset of $G$-orbit Closures}\label{S:Gorbs}

	Let $L:X$ be an algebraic group action. 
	By $I(L:X)$ we will denote the inclusion poset of $L$-orbit closures in $X$.

	\begin{Theorem}\label{T:GorbitinterectsY}
		Let $p: X\to G/H$ be a homogeneous fiber bundle over $G/H$. 
		Let $Y$ denote the fiber of $p$ at $eH$, that is, $Y:=p^{-1}(eH)$.
		If $G$ has an open orbit in $X$, then the posets $I(G:X)$ and $I(H:Y)$ are isomorphic. 
	\end{Theorem}
	
	\begin{proof}
		Let $O$ be an element from $I(G:X)$. 
		Then $O$ is of the form $\overline{G\cdot x}$ for some point $x$ in $X$. 
		We know from Lemma~\ref{L:GorbitinterectsY1} that $O \cap Y\neq \emptyset$. 
		We will show that $O\cap Y$ is actually an $H$-orbit closure. 
		Let $y$ be an element from $(G\cdot x)\cap Y$. 
		Clearly, we have $H\cdot y \subseteq (G\cdot x)\cap Y$. 
		We claim that this inclusion is actually equality, that is, $H\cdot y = (G\cdot x)\cap Y$. 
		Let $z$ be an element from $(G\cdot x)\cap Y$. 
		Then $z= g\cdot x$ for some $g\in G$. 
		Since both $z$ and $y$ are elements of the $G$-orbit $G\cdot x$, there exists $d\in G$ such that 
		$d\cdot z = y$.  
		Then by Lemma~\ref{L:GorbitinterectsY2}, we have $d\in H$. 
		In particular, we see that $z\in H\cdot y$ in $Y$. 
		This argument finishes not only the proof of our claim but also shows that the following map is well-defined:
		\begin{align*}
			\varphi : I(G:X) &\longrightarrow I(H:Y) \\  
			\overline{G\cdot x} &\longmapsto \overline{Y\cap (G\cdot x)}.
		\end{align*} 
		Let $K$ be an element of $I(H:Y)$. 
		Then $K= \overline{ H\cdot y}$ for some $y\in Y$. 
		Clearly, the $G$-orbit $G\cdot y$ intersects $Y$ along $K$. 
		Therefore, $\varphi$ is surjective. 
		Now let $K_1:= \overline{H\cdot y_1}$ and $K_2= \overline{H\cdot y_2}$ be two distinct $H$-orbit closures in $Y$.
		We claim that $(G\cdot y_1) \cap (G\cdot y_2) = \emptyset $. 
		Indeed, if these two $G$-orbits coincide, then since $y_i \in H$ ($i\in \{1,2\}$), by Lemma~\ref{L:GorbitinterectsY2}, there exists $d\in H$ 
		such that $y_1 = d\cdot y_2$. But this contradicts with our assumption that $K_1 \neq K_2$. 
		This proves our claim. Hence, $\varphi$ is injective.
		Thus we proved that $\varphi$ is a bijection between $I(G:X)$ and $I(H:Y)$. 
		We are now ready to prove that it is an order isomorphism, that is, 
		\[
		O_1 \subseteq O_2 \iff \varphi(O_1) \subseteq \varphi(O_2)
		\]
		for every $O_1,O_2$ from $I(G:X)$. 
		We will use induction on $\dim X$. 
		
		If $\dim X =1$, then we have one of the following two possibilities: 
		\begin{enumerate}
			\item[(a)] $\dim G/H = 0$,
			\item[(b)] $\dim G/H =1$
		\end{enumerate} 
		In the former case, since $G$ is connected, we have $G=H$. Then $G/H$ is a point and $Y=X$. 
		Thus, $I(G:X)=I(H:Y) = \{ X \}$. 
		In the latter case, for dimension reasons, we have $X= G/H$.
		It follows that $Y= eH$.
		In other words, we have $I(G:X) = \{X\}$ and $I(H:Y) = \{ Y \}$.
		Since both posets have single element, they are isomorphic.

		We now assume that our claim holds for every homogeneous fiber bundle $Z\to G/H$ such that $\dim Z = n$.
		Let $X$ be a homogeneous fiber bundle over $G/H$ such that $\dim X = n+1$. 
		Let $O_1$ and $O_2$ be two $G$-orbit closures from $I(G:X)$. 
		
		First, we assume that $O_1\subseteq O_2$. 
		Let us assume also that $\dim O_2 < n+1$.
		Notice that the restriction of $p$ to $O_2$ gives a surjective $G$-equivariant morphism $p|_{O_2} : O_2\to G/H$.
		Hence, $O_2$ stands as a homogeneous fiber bundle itself. 
		Since its dimension does not exceed $n$, by applying our inductive assumption, we readily see that the inclusion $\varphi(O_1)\subseteq \varphi(O_2)$ holds in $I(G:O_2)$. But this is a subposet of $I(G:X)$. Hence, the inclusion $\varphi(O_1)\subseteq \varphi(O_2)$ holds in $I(G:X)$ as well.
		If $\dim O_2 = n+1$, then the $G$-orbit $G\cdot x$ such that $\overline{G\cdot x} = O_2$ is the open orbit in $X$. 
		Hence, the intersection $(G\cdot x) \cap Y$ is open in $Y$. 
		In particular, we have $\varphi(O_2) = Y$. 
		Therefore, the inclusion $\varphi(O_1) \subseteq \varphi (O_2)$ holds in $I(G:X)$ in this case, also.
		
		Conversely, we assume that $\varphi(O_1) \subseteq \varphi (O_2)$. 
		We want to show that $O_1\subseteq O_2$. 
		To this end, we notice that the inclusion $\varphi(O_1) \subseteq \varphi (O_2)$ implies 
		that $G\cdot \varphi(O_1) \subseteq G\cdot \varphi (O_2)$. 
		Since we have the equality $\overline{G\cdot \varphi(O_i)} = O_i$ for $i\in \{1,2\}$, 
		our claim follows.
		Hence, the proof is complete. 
	\qed\end{proof}

\begin{Remark}
Thanks to Corollary~\ref{C:Gorbs}, in Theorem~\ref{T:GorbitinterectsY}, we may replace the hypothesis that $G$ has an open orbit by the following two assumptions: 1) $Y$ is irreducible, and 2) $H$ has an open orbit in $Y$.
Alternatively, we may use the following two assumptions: 1') $H$ is connected, and 2') $H$ has a dense orbit in $Y$.
To see that 1' and 2' together imply 1 and 2, we first notice that since $H$ is connected, the closure of an orbit of $H$ is irreducible. 
Since $H$ has a dense orbit $O\subseteq Y$, its Zariski closure, that is, $Y$ is irreducible.   
Secondly, since orbits of maximal dimension form an open set, we see that $O$ is open.
Hence, we showed that the assumptions 1' and 2' imply the assumptions 1 and 2.
\end{Remark}

\begin{Example}
Let $G$ be a connected reductive group.
Let $B$ be a Borel subgroup.
Let $\chi$ be a nontrivial character of $B$. 
As we mentioned before, the projection $G\times_B k_\chi \to G/B$ is a homogeneous (line) bundle over $G/B$.  
The fiber of $G\times_B k_\chi \to G/B$ at $eB$ is given by $\{e\} \times_B k_\chi$.
This is naturally isomorphic to the vector space $k$ on which $B$ acts via $\chi$. 
In characteristic 0, this action has only two orbits; $\{ 0 \}$ and $k\setminus \{0\}$. 
The former orbit is contained in the closure of the latter orbit. 
It follows that $G$ has two orbits on $G\times_B k_\chi$ with the same inclusion relationship. 
\end{Example}

\section{$G$-BSDH Varieties}\label{S:GBSDH}

	We begin with fixing some additional notation. 
	Hereafter, $G$ will denote a connected reductive group. 
	As usual, $B$ will denote a Borel subgroup of $G$, and $T$ will denote a maximal torus of $B$. 
	The Weyl group of $(G,T)$ is denoted by $W$. 
	The set of Coxeter generators of $W$ (relative to $B$) is denoted by $S$ (or by $S(G)$ when we need to make a distinction). 
	We let $\leq$ denote the Bruhat-Chevalley order. 
	The length function on $W$ will be denoted by $\ell$. 
	\medskip

	The unipotent radical of an algebraic group $H$ will be denoted by $R_u(H)$. 
	If $H$ is the Borel subgroup $B$, then we will use the letter $U$ instead of $R_u(H)$.
	The unique Borel subgroup that is opposite to $B$ is denoted by $B^-$.
	Also, its unipotent radical will be denoted by $U^-$. 
	\medskip
	
	Let $P$ be a parabolic subgroup of $G$. 
	If $P$ is {\em standard with respect to $B$}, that is to say $B\subseteq P$, then we write $P_J$ instead of $P$,
	where, $J\subseteq S$ is the set of Coxeter generators such that $B$ and the representatives of the elements of $J$ in $G$ generate $P_J$. 
	In this notation, the Weyl group of $P_J$ is denoted by $W_J$.
	The corresponding set of minimal length left coset representatives in $W$ is denoted by $W^J$. 
	For $s\in S$, the minimal parabolic subgroup generated by $B$ and the representative $n_s\in G$ of $s$ is denoted by $P_s$. 
	Then $P_s$ is given by the union $B\cup B s B$. 
	Notice that, in this union, we used $s$ instead of $n_s$.
	We will continue to follow this convention in the sequel to simplify our notation.
	\medskip

	A {\em word in $S$} is a finite sequence of not necessarily distinct elements from $S$. 
	If $\underline{w}$ denotes the word $(s_{i_1},\dots, s_{i_m})$ in $S$, then the product $w:=s_{i_1}\cdots s_{i_m}$ is an element of $W$. 
	In this notation, the corresponding Schubert variety is given by 
	\[
	X_{wB} = P_{s_{i_1}}P_{s_{i_2}}\cdots P_{s_{i_m}} / B.
	\]
	We say that a word $\underline{w}=(s_{i_1},\dots, s_{i_m})$ is a {\em reduced word} if its length $m$ equals the dimension of the corresponding Schubert variety, $m=\dim X_{wB}$. 
	\medskip

	We are now ready to discuss our generalized Schubert varieties. 
	We consider the morphism 
	\begin{align}\label{A:firstformofxi}
		\xi : G\times_B G/B &\longrightarrow G/B\times G/B \\ 
		[g,g'B] &\longmapsto (gB,gg'B). \notag
	\end{align}
	It is well-known that every closed irreducible $G$-stable subvariety of $G/B\times G/B$ is of the form $\xi(G\times_B X_{wB})$ for some $w\in W$,
	see~\cite[page 69]{BrionKumar}. 
	For this reason, a homogeneous fiber bundle of the form $G\times_B X_{wB}$, where $w\in W$, is called a {\em $G$-Schubert variety}.

	\medskip
	
	We recall the definition of the BSDH-varieties to set up our notation.
		Let $\underline{w}$ be a word from $S$. 
		If $\underline{w} = (s_{i_1},\dots, s_{i_m})$, then we have the product of minimal parabolic subgroups $\prod_{j=1}^m P_{i_j}$,
		where $P_{i_j} = B\cup B s_{i_j} B$ for $j\in \{1,\dots, m\}$.
		Then the BSDH-variety $X_{\underline{w}}$ is defined as the quotient variety $\prod_{j=1}^m P_{i_j} / B^m$, 
		where the (right) action of $B^m$ is given by 
		\[
		(p_1,\dots, p_m) \cdot (b_1,\dots, b_m) = (p_1 b_1, b_1^{-1} p_2 b_2 , \dots, b_{m-1}^{-1} p_m b_m),
		\]
		for $(b_1,\dots, b_m)\in B^m$ and $(p_1,\dots, p_m)\in P_{i_1}\times \cdots \times P_{i_m}$. 
		The image of an element $(p_1,\dots, p_m)$ of $P_{i_1}\times \cdots \times P_{i_m}$ 
		under the quotient map $\pi_{\underline{w}} : \prod_{j=1}^m P_{i_j} \to X_{\underline{w}}$
		will be denoted by $[p_1,\dots, p_m]$.
		\medskip
	The BSDH-variety $X_{\underline{w}}$ is closely related to the corresponding Schubert variety $X_{wB}$. 
	Indeed, the Schubert variety $X_{wB}$ is image of the natural product map,   
	\begin{align*}
		\mathbf{m} :\ X_{\underline{w}} &\longrightarrow G/B \\
		[p_1,\dots, p_m] & \longmapsto p_1\cdots p_mB. 
	\end{align*}
	It is well-known that~\cite[Theorem 3.4.3]{BrionKumar} if $\underline{w}$ is a reduced word, then $\mathbf{m}$ is a resolution of singularities of $X_{wB}$. 
	We note in passing that in~\cite{BrionKumar} a BSDH-variety associated with $\underline{w}$ is denoted by $Z_{\underline{w}}$.

\begin{Definition}
A {\em $G$-BSDH variety} is a homogeneous fiber bundle of the form $G\times_B X_{\underline{w}}$, where $\underline{w}$ is a word in $S$. 
\end{Definition}

We are now ready to prove our Theorem~\ref{T:toroidal}.
Let us recall its statement for convenience. 
\medskip
	
Let $X:=G\times_B X_{\underline{w}}$ be a $G$-BSDH variety. 
If a $B$-stable divisor $D$ in $G\times_B X_{\underline{w}}$ contains a $G$-orbit, then $D$ is $G$-stable. 
\medskip

\begin{proof}[Proof of Theorem~\ref{T:toroidal}]
We fix a $B$-linearized very ample line bundle $L\to X_{\underline{w}}$. 
Let $Y\subseteq \mathbb{A}^{N+1}$ denote the affine cone for the corresponding projective embedding $X_{\underline{w}}\to \mathbb{P}^N$. 
Let $Y_0:= Y\setminus \{0\}$. 
Then $G\times Y_0$ is a quasi-affine variety. 
Let $q:Y_0\to X_{\underline{w}}$ denote the quotient map. 
This is a smooth morphism. 
We consider the following morphisms:
\begin{equation*}
\begin{aligned}[t]
  p:G\times Y_0 &\to G\times X_{\underline{w}}\\
  (g,x)              &\mapsto (g,q(x)),
\end{aligned}
\qquad\text{and}\qquad
\begin{aligned}[t]
  \pi: G\times X_{\underline{w}} &\to G\times_B X_{\underline{w}}\\
  (g,q(x))                &\mapsto [(g,q(x))].
\end{aligned}
\end{equation*}
Let $\tilde{p}: G\times Y_0 \to G\times_B X_{\underline{w}}$ denote their composition as in Figure~\ref{F:composition}.
\begin{figure}[htp]
\begin{center}
\begin{tikzcd}[row sep = huge, column sep = huge,every label/.append  style={font=\normalsize}]
\makebox{$G\times Y_0$} \arrow[r, "p",  rightarrow] \arrow[rd, "\tilde{p}", rightarrow, dashed]  & \makebox{$G\times X_{\underline{w}}$}  \arrow[d, "\pi", rightarrow] \\
   & \makebox{$G\times_B X_{\underline{w}}$}  \arrow[d, "\alpha_1",  rightarrow]\\
   & \makebox{$G/B$}
\end{tikzcd}
\end{center}
\caption{A diagram of $G$-equivariant quotient maps.}
\label{F:composition}
\end{figure}
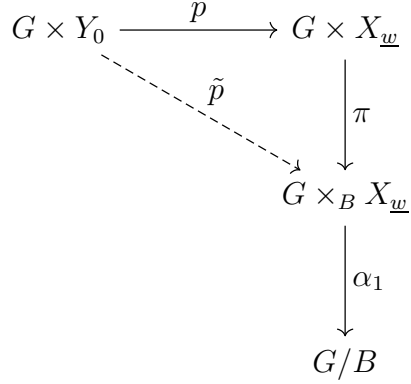 
It is easy to check that both $p$ and $\pi$ are smooth, $G$-equivariant quotient morphisms. 
It follows that $\tilde{p}$ is a smooth, $G$-equivariant, and surjective morphism of quasi-projective varieties as well.
It is also evident that $\tilde{p}(\{e\}\times Y_0) = \{e\} \times_B X_{\underline{w}}$. 
\medskip	

Now let $D$ be a $B$-stable divisor in $G\times_B X_{\underline{w}}$. 
Then, the preimage $\tilde{p}^{-1}(D)$ is a $B$-stable divisor in $G\times Y_0$.
Let $O$ be a $G$-orbit contained in $D$. 
Then there exists a point $[h,x]\in D$, where $h\in G$ and $x\in X_{\underline{w}}$, such that $O=G\cdot [h,x]$.
Since $G$ acts on the first coordinate, we have $O=G\times_B \{x\}$. 
Notice that $\pi (G\times \{x\}) = G\times_B \{x\}=O$. 
Hence, we have the following containments of quasi-affine varieties:
\[
G\times \{ x\} \subset \tilde{p}^{-1} (O) \subset \tilde{p}^{-1} (D).
\] 
It follows that the coordinate ring $k[G]$ is a quotient of the coordinate ring of $\tilde{p}^{-1}(D)$. 
\medskip

We proceed to make a related observation. 
We will consider the first projection,  
\begin{align*}
\alpha_1: G\times_B X_{\underline{w}}&\rightarrow G/B \\
[g,x]&\mapsto gB,
\end{align*}
which is a $G$-equivariant morphism. 
Then we have $\alpha_1^{-1}(eB)=\{e\}\times_B X_{\underline{w}}$.
Let us call this preimage the {\em special fiber}.
We know that every $G$-orbit in $G\times_B X_{\underline{w}}$ intersects the special fiber.
Let $F:=D\cap (\{e\}\times_B X_{\underline{w}})$. 
Since $D$ contains a $G$-orbit, $F$ is nonempty. 
At the same time, since $D$ is a divisor in $G\times_B X_{\underline{w}}$, exactly one of the following two cases may occur:
\begin{enumerate}
\item[(1)] $F$ is a divisor in $\{e\}\times_B X_{\underline{w}}$, or 
\item[(2)] $\{e\}\times_B X_{\underline{w}}$ is contained in $D$, so $F = \{e\} \times_B X_{\underline{w}}$.
\end{enumerate}

We will show that (2) is not possible. Towards a contradiction, let us assume that (2) holds. 
Then, we see that $\{e\}\times Y_0$ is a subvariety of $\tilde{p}^{-1}(D)$.
This means that the coordinate ring $k[Y_0]$ is a quotient of the coordinate ring of $\tilde{p}^{-1}(D)$.
We know the following facts: 1) $(e,x) \in \{e\}\times Y_0$, 2) $D$ contains the $G$-orbit of $(e,x)$, 3) $\tilde{p}$ is $G$-equivariant,
and 4) $G$ is connected. Therefore, the irreducible component of $\tilde{p}^{-1}(D)$ that contains $\{e\} \times Y_0$ also contains the $G$-orbit of $(e,x)$.
But this implies that $G\times \{x\}$ is a subvariety of the irreducible component of $\tilde{p}^{-1}(D)$ that contains $\{e\}\times Y_0$. 
Now we know that $\tilde{p}^{-1}(D)$ contains two closed subsets $Z_1:=G\times \{x\}$ and $Z_2:=\{e\}\times Y_0$ which intersect along the point $(e,x)$. 
It follows that $\tilde{p}^{-1}(D)$ contains a subvariety that is isomorphic to 
\begin{align*}
J:= (G\times \{x\}) \ \times \  (\{e\}\times Y_0).
\end{align*}
Since the coordinate ring of $J$ is isomorphic to $k[G]\otimes k[Y_0]$, we see that the coordinate ring $k[\tilde{p}^{-1}(D)]$ contains a subring that is isomorphic to $k[G]\otimes k[Y_0]$. 
Nonetheless, this ring is isomorphic to the coordinate ring of the ambient variety $G\times Y_0 \supsetneq \tilde{p}^{-1}(D)$. 
Hence, we obtain a contradiction from the fact that $\tilde{p}^{-1}(D)$ is only a divisor in $G\times Y_0$. 
\medskip

We proceed with (1).
Since both $D$ and $\{e\}\times_B X_{\underline{w}}$ are $B$-stable, we see that $F:= D\ \cap\ (\{e\}\times_B X_{\underline{w}})$ is a $B$-stable divisor in $\{e\}\times_B X_{\underline{w}}$.
Let $F_0$ denote the image of $F$ under the second projection $\alpha_2: G\times_B X_{\underline{w}}\to X_{\underline{w}}$. 
Let $O$ denote the $G$-orbit of some point $[e,z]$ in $F$. 
Let $\overline{O}$ denote the Zariski closure of $O$ in $X$.
It is easy to that $\pi^{-1}(\overline{O})$ is the divisor $G\times F_0$ in $G\times X_{\underline{w}}$.
It follows that $\pi^{-1}(D)$ is contained in $G\times F_0= (\pi \circ \alpha_2)^{-1}(F_0)$. 
But $\pi : G\times X_{\underline{w}}\to G\times_B X_{\underline{w}}$ is a smooth morphism, and $D$ is a prime divisor. 
Hence, $\pi^{-1}(D)$ is a prime divisor. 
Clearly, $G\times F_0$ is a prime divisor as well. 
It follows that we have the equality, 
\begin{align*}
\pi^{-1}(D) = G\times F_0.
\end{align*}
In particular, we see that $\pi^{-1}(D)$ is $G$-stable. 
Since $\pi$ is $G$-equivariant and $D= \pi(\pi^{-1}(D))$, the proof of our theorem is finished.
\qed\end{proof}

\begin{Corollary}\label{C:main1toroidal}
Let $X_{wB}$ be a Schubert variety in $G/B$.
Let $X$ denote the corresponding $G$-Schubert variety $G\times_B X_{wB}$. 
Let $D'$ be a $B$-stable divisor in $X$. 
If $D'$ contains a $G$-orbit, then $D'$ is $G$-stable. 
\end{Corollary}
	
\begin{proof}
We fix a reduced word $\underline{w}$ of $w$. 
Since the map $\mathbf{m} : X_{\underline{w}}\to X_{wB}$ is a resolution of singularities, the induced map 
$1\times \mathbf{m} : G\times_B X_{\underline{w}} \to X$ is a resolution of singularities as well. 
At the same time, $1\times \mathbf{m}$ is a $G$-equivariant morphism, where the action of $G$ is given by 
\[
g \cdot [h,x] = [gh,x], \qquad \text{where $g\in G$},
\]
and $[h ,x]$ represents either an element of $G\times_B X_{\underline{w}}$ or an element of $G\times_B X_{wB}$.
Now, since $1\times \mathbf{m}$ is $B$-equivariant and birational, we see that the preimage of every $B$-stable divisor $D'$ in $X$ is a $B$-stable divisor $D$ in $G\times_B X_{\underline{w}}$. 
Let $D'$ be a $B$-stable divisor in $X$. 
We claim that if $D'$ contains a $G$-orbit, then the divisor $D:=(1\times \mathbf{m})^{-1}(D')$ contains a $G$-orbit as well.  
To see this, let $y\in D'$ be a point such that $G\cdot y \in D'$. 
Let $z\in (1\times \mathbf{m})^{-1}(y)$. 
Assume towards a contradiction that there exists $t:= g\cdot z \in G\cdot z$ such that $t\notin D$.
Then we have $(1\times \mathbf{m}) (t) \notin D'$.
But since $1\times \mathbf{m}$ is $G$-equivariant, we obtain the following contradictory membership:
\[
(1\times \mathbf{m}) (t) = g\cdot (1\times \mathbf{m})(z) = g\cdot y \in D'.
\]
This contradiction shows $D$ contains a $G$-orbit as well. 
Then, it follows from Theorem~\ref{T:toroidal} that $D$ is $G$-stable. 
By using the $G$-equivariance of $1\times \mathbf{m}$ once more, we see that $D'$ is $G$-stable as well. 
This finishes the proof of our assertion.
\qed\end{proof}

\begin{Corollary}\label{C:main1regular}
Let $\underline{w}$ be a word in $S$. 
Let $w$ denote the corresponding element in $W$. 
If $G\times_B X_{\underline{w}}$ is a spherical variety, then both of the varieties $G\times_B X_{\underline{w}}$ and $G\times_B X_{wB}$ are regular $G$-varieties. 
\end{Corollary}
	
\begin{proof}
Let $X\in \{ G\times_B X_{\underline{w}},G\times_B X_{wB}\}$. 
By Theorem~\ref{T:toroidal} and its Corollary~\ref{C:main1toroidal}, we know that for every $G$-orbit $Y$ in $X$, we have $\mc{D}_Y =\emptyset$. 
Here, $\mc{D}_Y$ is the notation introduced in (\ref{A:notationtoroidal}). 
Hence, we have $\mc{D}_{all} = \emptyset$. 
In other words, $X$ is a toroidal variety.
The rest of the proof follows from Theorem~\ref{T:BienBrion}. 
\qed\end{proof}

\section{Toroidal Homogeneous Fiber Bundles}\label{S:Toroidal}

	In light of our Corollary~\ref{C:main1regular}, It is important to find a characterization of the spherical $G$-BSDH varieties. 
	Finding obstructions that prevent variety $G\times_B X_{\underline{w}}$ from being a spherical $G$-variety is not difficult;
	dimension provides the first obstruction. 
	
	\begin{Proposition}
		Let $\underline{w}$ be a word in $S$, and let $w$ be the corresponding element of $W$. 
		If $\ell(w) > |S|$, then $G\times_B X_{\underline{w}}$ is not spherical. 
	\end{Proposition}
	
	\begin{proof}
	First, we note that $\dim T = |S|$ and $\dim X_{\underline{w}} \geq \dim X_{wB} = \ell(w)$. 
		Let $\tilde{z}_0$ be a point from $X_{\underline{w}}$. 
		Since $\ell(w) > |S|$, the dimension of $X_{wB}$, hence the dimension of $X_{\underline{w}}$ is at least $\dim T+1$. 
		Let us compute the dimension of the homogeneous fiber bundle $G\times_B X_{\underline{w}}$, 
		\begin{align}\label{A:dimensioncalculation}
			\dim G\times_B X_{\underline{w}} &= \dim G + \dim X_{\underline{w}} - \dim B \notag \\ 
			&= \dim U^- + \dim X_{\underline{w}}.
		\end{align}
		It follows that we have the inequality $\dim G\times_B X_{\underline{w}} > \dim U^- + \dim T = \dim B$.
		Since the dimension of a $B$-orbit is at most $\dim B$, we see that $\dim G\times_B X_{\underline{w}}$ cannot have an open $B$-orbit. 
		This finishes the proof of our assertion. 
	\qed\end{proof}

	We continue with a finer analysis of the number of parameters on which a family of $G$-orbits in $G\times_B X_{\underline{w}}$ may depend.

\begin{Proposition}\label{P:modalityofGBSDH}
Let $\underline{w}$ be a word in $S$. Then we have 
\begin{align*}
\text{mod} (G:G\times_B X_{\underline{w}}) = \text{mod} (B:X_{\underline{w}}).
\end{align*}
In particular, if $\underline{w}$ is a reduced word of length $l$, then we have $\text{mod} (B:X_{\underline{v}})=0$ for every subword $\underline{v}$ of length $l-1$ if and only if we have $\text{mod} (G : G\times_B X_{\underline{w}}) = 0$.
\end{Proposition}
	
\begin{proof}
Recall that the modality of an action $G:X$ is the number $\max_{Y\subseteq X} \ \text{tr.deg } k(Y)^G$,
where the maximum is taken on the set of all $G$-stable irreducible subvarieties $Y$ of $X$. 
Recall also that every $G$-orbit in $G\times_B X_{\underline{w}}$ is of the form $G\times_B Z$, where $Z$ is a $B$-orbit in $X_{\underline{w}}$. 
Since $k(G\times_B Z) \cong k(G\times Z)^B$, and since the left action of $G$ and the right action of $B$ commute with each other, 
we see that 
\[
{}^Gk(G\times_B Z)= {}^G (k(G\times Z)^B)= ({}^G k(G\times Z))^B= k(Z)^B.
\]
Thus, the transcendence degree of ${}^Gk(G\times_B Z)$ is equal to the transcendence degree of $k(Z)^B$. 
In light of the correspondence between the $G$-orbits in $G\times_B X_{\underline{w}}$ and the $B$-orbits in $X_{\underline{w}}$,
we see that 
$$
\max_{Y:\text{ $B$-stable irr. subv. of $X_{\underline{w}}$}} \text{tr.deg } k(Y)^B 
= \max_{Y:\text{ $G$-stable irr. subv. of $G\times_B X_{\underline{w}}$}}  \text{tr.deg } {}^Gk(G\times_B Y).
$$ 
This finishes the proof of our first assertion.

To prove our second assertion we note the fact that if $\underline{w}$ is a reduced word, then the $B$-equivariant surjective proper map, 
$\mathbf{m}:X_{\underline{w}}\to X_{wB}$ is a birational morphism. 
Let $j\in\{1,\dots, l\}$. 
If $\underline{w}$ is given by $\underline{w}:=(s_{i_1},\dots, s_{i_l})$, then let $\underline{w}^j$ denote word defined by suppressing the $j$-th entry,
\[
\underline{w}^j := (s_{i_1},\dots, s_{i_{j-1}},s_{i_{j+1}},\dots, s_{i_l}). 
\]
By~\cite[Proposition 2.2.6 (ii)]{Brion_Lectures}, we see that $X_{\underline{w}}$ is given by the union
\begin{align}\label{A:therighthandsideismodality0}
X_{\underline{w}} = \mathbf{m}^{-1}(C_{wB}) \cup \bigcup_{j=1}^{l} X_{\underline{w}^j},
\end{align}
where $C_{wB}$ is the open $B$-orbit in the Schubert variety $X_{wB}$. 
At the same time, we know from~\cite[Proposition 2.2.1 (iv)]{Brion_Lectures} that $\mathbf{m}^{-1}(C_{wB})$ is isomorphic to $C_{wB}$ via the restriction of the $B$-equivariant morphism, $\mathbf{m}$. 
Since $\text{mod} (B:X_{\underline{w}^j})=0$ if and only if $X_{\underline{w}^j}$ has only finitely many $B$-orbits, we see from (\ref{A:therighthandsideismodality0}) that $X_{\underline{w}}$ has only finitely many $B$-orbits. 
Clearly, this argument is reversible. 
Now, the rest of the proof follows from our first assertion. 
This finishes the proof of our proposition.
\qed\end{proof}

\begin{Example}
Let $G$ denote $GL(3,\C)$. 
Let $B$ denote the Borel subgroup of upper triangular matrices in $G$.
Let $T$ denote the maximal diagonal torus in $B$. 
The Weyl group of $(G,T)$ is the symmetric group $\mathbf{S}_3$. 
The set of simple reflections is given by $S=\{s_1,s_2\}$, where $s_i$ (for $i\in \{1,2\}$) is the simple transposition that interchanges $i$ and $i+1$. 
Let $w_0$ denote $s_1s_2s_1$. 
We will consider the reduced word $\underline{w_0}:=(s_1,s_2,s_1)$.
For $i\in \{1,2\}$, let $P_{s_i}$ denote the standard parabolic subgroup of $GL(3,\C)$ that is generated by $B$ and $s_i$. 
Now, the three subwords of $\underline{w_0}$ of length 2 are $(s_1,s_2)$, $(s_2,s_1)$, and $(s_1,s_1)$. 
Then the corresponding BSDH-varieties are given by $P_{s_1}\times_B P_{s_2}/B$, $P_{s_2}\times_B P_{s_1}/B$, 
and $P_{s_1}\times_{B} P_{s_{1}}/B\simeq P_{s_1}/B\times P_{s_1}/B$, respectively. 
The first two BSDH-varieties are isomorphic to the Hirzebruch surface $\mathbb{P}(\mathcal{O}\oplus \mathcal{O}(-1))$.
Here, $\mathcal{O}$ is the structure sheaf of $\mathbb{P}^1$. 
Let $\underline{v}:=(s_1,s_1)$. 
Then $X_{\underline{v}}$ is isomorphic to $\mathbb{P}^1\times \mathbb{P}^1$ as an $L_{\{s_1\}}$-variety, where $L_{\{s_1\}}$ 
is the standard Levi subgroup of $P_{\{s_1\}}$. The maximal torus $T\subset L_{\{s_1\}}$ has an open orbit in $\mathbb{P}^1\times \mathbb{P}^1$.  
We conclude from this calculation that for every (not necessarily reduced) subword $\underline{v}$ of length $2$ of $\underline{w_0}=(s_1,s_2,s_1)$, 
the corresponding BSDH-variety $X_{\underline{v}}$ contains only finitely many $B$-orbits. 
\end{Example}

We know from Karuppuchamy's work (\cite[Theorems 2 and 4]{Karuppuchamy}) that a Schubert variety $X_{wB}$ is a toric variety if and only if $w$ is a product of distinct simple reflections from $S$. 
It turns out that $X_{wB}$ being a toric variety is equivalent to the $G$-Schubert variety $G\times_B X_{wB}$ being a spherical variety. 
To state this more precisely, we will briefly review some results which are originally due to Luna. 
\medskip

	\begin{Definition}
	Let $P$ be a parabolic subgroup. 
	Let $L$ be a Levi factor of $P$.
	We denote by $B_L$ a Borel subgroup of $L$.  
	A normal $P$-variety $Z$ is called a {\em spherical $P$-variety} if the $Z$ is a spherical $L$-variety. 
	A spherical $P$-variety $Z$ is said to be {\em wonderful} if it possesses the following properties:
	\begin{enumerate}
		\item[(1)] $Z$ is smooth and complete; 
		\item[(2)] $Z$ contains exactly one closed $P$-orbit $Z_0$;
		\item[(3)] every irreducible $B_L$-stable closed subvariety $Z'\subseteq Z$ containing $Z_0$ is actually $P$-stable. 
	\end{enumerate}
	\end{Definition}
	Notice that the properties (1)--(3) are equivalent to saying that $Z$ is smooth, complete, simple, and toroidal. 
	
	\begin{Lemma}\cite[Proposition 6.3]{Avdeev2015}\label{L:Avdeev's2}
		Let $Z$ be a $P$-variety and consider the $G$-variety $X= G\times_P Z$. 
		Then we have 
		\begin{enumerate}
			\item $Z$ is a spherical $P$-variety if and only if $X$ is a spherical $G$-variety.
			\item $Z$ is a wonderful $P$-variety if and only if $X$ is a wonderful $G$-variety.
		\end{enumerate}
	\end{Lemma}

	We are ready to prove the main result of this section.
	
	\begin{Proposition}\label{P:fourequivalent}
		Let $\underline{w}$ be a reduced word in $S$. 
		Then the following statements are equivalent:
		\begin{enumerate}
			\item $X_{\underline{w}}$ is a toric variety;
			\item $X_{wB}$ is a toric variety; 
			\item $G\times_B X_{\underline{w}}$ is a spherical $G$-variety; 
			\item $G\times_B X_{wB}$ is a spherical $G$-variety. 
		\end{enumerate}
	\end{Proposition}
	
	\begin{proof}
		The two sided implications (1)$\Leftrightarrow$(3) and (2)$\Leftrightarrow$(4) follow from Lemma~\ref{L:Avdeev's2} where $P$ is $B$. 
		Since $\underline{w}$ is a reduced word, the canonical $B$-equivariant map $X_{\underline{w}}\to X_{wB}$ is a birational morphism. 
		Hence, $X_{\underline{w}}$ is a toric variety if and only if $X_{wB}$ is a toric variety. 
		This finishes the proof of our assertion.  
	\qed\end{proof}

	We are now ready to remove the reducedness assumption on $\underline{w}$. 
	Recall that our Theorem~\ref{T:main2} is the following statement: 
	\medskip

Let $\underline{w}$ be a word in $S$. Then $G\times_B X_{\underline{w}}$ is a wonderful variety if and only if $X_{\underline{w}}$ is a toric variety. 
Furthermore, if $\underline{w}$ is a reduced word in $S$, then $G\times_B X_{\underline{w}}$ is a wonderful variety if and only if $X_{wB}$ is a toric variety. 
	
\medskip

\begin{proof}[Proof of Theorem~\ref{T:main2}.]
By~\cite[Proposition 2.2.1]{BienBrion}, we know that a smooth and complete $G$-variety is regular if and only if it is spherical without color. 
It follows from this fact that a regular $G$-variety is spherical if it is complete. 
Clearly, every $G$-BSDH variety is complete.
We mentioned in Subsection~\ref{SS:toroidalregular} also that a regular $G$-variety $X$ is a wonderful variety if and only if $X$ has a unique closed $G$-orbit. 
Therefore, to prove our first claim it suffices to show that $G\times_B X_{\underline{w}}$ is a regular variety with unique closed orbit if and only if the BSDH-variety $X_{\underline{w}}$ is a toric variety. 
		
By Lemma~\ref{L:Avdeev's2} we know that $X_{\underline{w}}$ is a toric variety if and only if $G\times_B X_{\underline{w}}$ is spherical. 
Since $G\times_B X_{\underline{w}}$ is smooth, by Corollary~\ref{C:main1regular}, $G\times_B X_{\underline{w}}$ is regular variety. 
To show that $G\times_B X_{\underline{w}}$ is a simple $G$-variety, we will use Theorem~\ref{T:GorbitinterectsY}. 
Since $G\times_B X_{\underline{w}}$ is spherical, $G$ has an open orbit in $G\times_B X_{\underline{w}}$. 
Thus, by Theorem~\ref{T:GorbitinterectsY}, the poset of $G$-orbit closures in $G\times_B X_{\underline{w}}$ is isomorphic to the poset of $B$-orbit closures in $X_{\underline{w}}$. Since there is a unique $B$-fixed point in a BSDH-variety, we see that $G$ has a unique closed $G$-orbit in $G\times_B X_{\underline{w}}$. This finishes the proof of our first assertion. 
		
The proof of our second assertion follows from the first part combined with Proposition~\ref{P:fourequivalent}. 
\qed\end{proof}

	Evidently, wonderful varieties have a special place among spherical varieties. 
	We conclude this section by a discussion of the ``ranks'' of our wonderful varieties. 
	\medskip

Let $X$ be a spherical $G$-variety with open orbit denoted by $X_0$. 
Let $H$ be the stabilizer subgroup ${\rm Stab}_G(x_0)$, where $x_0\in X_0$. 
Then we have $k(X) = k(G/H)$. 
Under the natural action of $B$, $k(G/H)$ decomposes into eigenspaces. 
We denote by $\MO(X)$ the lattice of the $B$-weights of the $B$-eigenvectors in $k(G/H)$.
This is a sublattice of the character lattice of $(G,T)$. 
The rank of $\MO(X)$ is called the {\em spherical rank of $X$}. 
There is a distinguished basis for the $\Q$-vector space $\MO(X)\otimes_\Z \Q$ consisting of certain primitive elements of $\MO(X)$. 
The elements of this basis are called the {\em spherical roots of $G/H$}.
It is observed in~\cite{Knop1996} that this basis generates $\MO(X)$ as a lattice if and only if $X$ is a wonderful variety. 
Let us now continue with the assumption that $X$ is a wonderful $G$-variety.
Then the {\em rank} of $X$ is defined as the number of irreducible components of the boundary, $X\setminus X_0$. 
Each such component is a $G$-stable prime divisor. 
In~\cite[Remark 3.17]{Avdeev2015}, Avdeev points out that the number of irreducible boundary divisors is equal to the spherical rank of $X$.

\begin{Proposition}
Let $\underline{w}$ be a reduced word in $S$. 
If $X_{\underline{w}}$ is a toric variety, then the poset of $G$-orbit closures in $G\times_B X_{\underline{w}}$ is isomorphic to the boolean lattice of all 
subsets of the set $\{1,\dots, \ell\}$, where $\ell = \dim X_{\underline{w}}$. 
In particular, the spherical rank of $G\times_B X_{\underline{w}}$ equals $\ell = \dim X_{\underline{w}}$.
\end{Proposition}
	
\begin{proof}
Since $\underline{w}$ is a reduced word, if $X_{\underline{w}}$ is a toric variety, then $\mathbf{m} : X_{\underline{w}}\to X_{wB}$ is an isomorphism.
In particular, we have $\dim X_{\underline{w}} = \ell = \dim X_{wB}$. 
In fact, since $X_{\underline{w}}$ is $B$-equivariantly isomorphic to $X_{wB}$, we see that the number of $B$-orbits in $X_{\underline{w}}$ is finite, implying $\text{mod}(B:X_{\underline{w}}) = 0$. 
It follows from Proposition~\ref{P:modalityofGBSDH} that $G\times_B X_{\underline{w}}$ has an open $G$-orbit. 
Now by Theorem~\ref{T:GorbitinterectsY}, the poset of $G$-orbit closures in $G\times_B X_{\underline{w}}$ is isomorphic to the poset of $B$-orbit closures in $X_{\underline{w}}$. 
Since we have $X_{\underline{w}} \cong X_{wB}$ by the $B$-equivariant map $\mathbf{m}$, we see that the poset of $B$-orbit closures in $X_{\underline{w}}$ is given by the lower interval $[e,w]$ of the Bruhat poset $(W,\leq)$. In particular, the number of $B$-stable prime divisors of $X_{\underline{w}}$ is equal to the co-atoms of the interval $[e,w]$. But $w$ is a product of $\ell$ distinct simple reflections. 
Hence, deleting a simple reflection from $w$ gives a coatom. 
It follows that the number of distinct simple reflection in $w$ is the number of $G$-stable prime divisors of $G\times_B X_{\underline{w}}$. 
This finishes the proof of our assertion. 
\qed\end{proof}

The arguments that we used in the proof of our previous proposition can be used for proving the following result.
We omit its details. 
	
\begin{Proposition}
Let $\underline{w}$ be a reduced word in $S$. 
Let $w$ denote the element associated with $\underline{w}$ in $W$. 
If $X_{\underline{w}}$ has finitely many $B$-orbits, then the poset of $G$-orbit closures in $G\times_B X_{\underline{w}}$ is isomorphic to the boolean lattice structure on the set of $T$-fixed points, $X_{\underline{w}}^T$.	
\end{Proposition}

\section{Parabolic Decompositions and the Spherical Schubert Varieties}\label{S:PD}
	
In this section $G$ denotes a connected semisimple algebraic group over an algebraically closed field $k$ of arbitrary characteristic. 
We proceed with establishing our notational conventions.
	
\begin{Notation}
Earlier, we used the letter $S$ to denote the set of Coxeter generators of $W$ determined by the pair $(B,T)$.
Hereafter, when confusion is unlikely, we will use $S$ to denote the set of simple roots corresponding to this set of Coxeter generators as well.
We will follow the standard convention that if $\alpha$ is a simple root, then the corresponding simple reflection is denoted $s_\alpha$. 
However, if the set of simple roots are given by an ordered set $\{\alpha_1,\dots, \alpha_n\}$, then the corresponding simple reflections will be written in the form $\{s_1,\dots, s_n\}$, where $s_i$ corresponds to $\alpha_i$ for $1\leq i \leq n$.
The root system of the pair $(G,T)$ will be denoted by $R$.  
We use the notation $R^+$ for the system of positive roots determined by $S$. 
If $\beta$ is a (positive) root, then its support, denoted by $supp(\beta)$, is the set of simple roots that appear in $\beta$ as a summand.
\end{Notation}
\medskip

Let $\beta$ be a root from $R$. 
The root subgroup associated with $\beta$ is denoted by $U_\beta$.  
It is given by the image of the isomorphism $x_{\beta}: \mathbb{G}_{a}\to U_{\beta}$ satisfying  
\[
tx_{\beta}(a)t^{-1}=x_{\beta}(\beta(t)a)\qquad\text{for $t\in T$ and $a\in \mathbb{G}_{a}$}.
\]
Let $w\in W$. Let $X_{wB}$ be the corresponding Schubert variety. 
The stabilizer subgroup $P(w)={\rm Stab}_{G}(X_{wB})$ is always a parabolic subgroup since it contains the Borel subgroup $B$. 
In particular, $P(w)$ is a {\em standard} parabolic subgroup of $G$. 
This means that there exists a subset $J(w)\subseteq S$ such that $P(w)=P_{J(w)}$, where $P_{J(w)}$ is generated by $B$ and $\{s_\alpha :\ \alpha \in J(w)\}$.	
Combinatorially speaking, the set $\{s_\alpha :\ \alpha \in J(w)\}$ is the {\em left descent set of $w$}, that is, the set of simple reflections $s_\alpha$ in $W$ such that $\ell(s_\alpha w)<\ell(w)$. 
The {\em standard Levi factor} of $P_{J(w)}$ is the unique Levi subgroup $L(w)$ such that $T\subseteq L(w)$.
Let $L(w)$ denote the standard Levi factor of $P_{J(w)}$. Then the intersection $B\cap L(w)$, denoted $B_{L(w)}$, is a Borel subgroup of $L(w)$. 
\medskip

	For $w\in W$, we define 
	\[
	R^{+}(w^{-1}):=\{\beta\in R^{+}: w^{-1}(\beta)\in R\setminus R^+\}.
	\] 
	Note that $J(w)=R^{+}(w^{-1})\cap S$.
	We now define two special subsets of $R^+(w^{-1})$,
	\begin{align*}
	R_{1}:=\mathbb{Z}J(w)\cap R^{+}\qquad\text{and}\qquad R_{2}:=\{\beta\in R^{+}(w^{-1}) : {\text supp}(\beta)\nsubseteq J(w)\}. 
	\end{align*}	
	Recall that $\# R^+(w^{-1})  = \ell(w)$. 
	Then the key observation of this section is that $R_1$ and $R_2$ give a partitioning of $R^+(w^{-1})$: 
	\[
	R^{+}(w^{-1})=R_{1}\sqcup R_{2}.
	\] 
	Let $R_{1}=\{\beta_{1},\beta_{2},\ldots, \beta_{k}\}$ and $R_{2}=\{\beta_{k+1},\ldots, \beta_{r}\}$, where $\ell(w)=r$.
	It will be important for our purposes to keep in mind that the cardinality of $R_1$ equals the dimension of the flag variety of $L(w)$.
	Indeed, the elements of $R_1$ correspond to the root subgroups of the unipotent radical of $B_{L(w)}$. 
	\medskip

	Now we will prove a special case of the first announced theorem of our paper.

	\begin{Theorem}\label{thm6.1}
	Let $w\in W$. 
	Then the associated Schubert variety $X_{wB}$ is a spherical $L(w)$-variety such that $\dim B_{L(w)}= \dim X_{wB}$ if and only if $w$ can be written as 
	\[
	w=w_{0,J(w)}c, 
	\]
	where $w_{0,J(w)}$ is the longest element of $W_{J(w)}$ and $c$ is a Coxeter element such that $\ell(w)=\ell(w_{0,J(w)})+\ell(c)$. 
	\end{Theorem}
	
\begin{proof}
$(\Rightarrow)$
Assume that $X_{wB}$ is an $L(w)$-spherical variety such that $\dim B_{L(w)}= \dim X_{wB}$.
Then there exists a point $\xi$ in $X_{wB}$ such that $\overline{B_{L(w)}\xi}=X_{wB}$ and $\dim({\rm Stab}_{B_{L(w)}}(\xi))=0$. 
In fact, since $B_{L(w)}$ is a subgroup of $B$, $\xi$ is an element of the open $U$-orbit in $X_{wB}$. 
Now notice that 
\begin{align}\label{A:dimofflagofL(w)}
\dim(X_{wB})=\# R_{1}+ \# R_{2}=\dim (B_{L(w)})=\dim(T) + \dim(U_{w}),
\end{align}
where $U_{w}$ denotes the unipotent radical of $B_{L(w)}$.  
Since the dimension of the flag variety of $L(w)$ is given by 
\begin{align*}
\dim(U_{w})=\ell(w_{0,J(w)})=\# R_{1}, 
\end{align*}
it follows from the equality in (\ref{A:dimofflagofL(w)}) that $\# R_{2}=\dim(T)$. 
In other words, we have $k=r-n$, where $n={\rm rank}(G)$.
\medskip

Let $U^{*}=\{u\in U: x_{\beta}(u)\neq 0\ \text{ for all }\beta \in R^{+}\}$. 
Since $U^{*}wB/B$ is open in $X_{wB}$, and since $X_{wB}$ is irreducible, the following intersection is nonempty: 
\begin{align*}
B_{L(w)}\xi\ \cap\ U^*wB/B\neq \emptyset.
\end{align*}
Note that $U^{*}wB/B=U^{w}wB/B$, where $U^{w}:=\prod_{j=1}^{r} U_{\beta_{j}}^{*}$ and $U_{\beta_{j}}^{*}=U_{\beta_{j}}\setminus \{1\}$.
We write $\xi$ accordingly, as follows: 
\begin{align*}
\xi=\prod_{j=1}^{r}x_{\beta_{j}}(a_{j})wB/B
\end{align*}
for some $a_{1},\ldots, a_{k}\in \mathbb{G}_{a}$ and $a_{k+1},\ldots, a_{r}\in \mathbb{G}_{a}\setminus \{0\}$. 
Then, without loss of generality, we may assume that  
\begin{align}\label{A:wlog}
\xi=\prod_{j=k+1}^{r}x_{\beta_{j}}(a_{j})wB/B\qquad \text{for some $a_{j}\in \mathbb{G}_{a}\setminus \{0\}$.}
\end{align} 
Indeed, if $\xi=\prod_{j=1}^{r}x_{\beta_{j}}(a_{j})wB/B$ for some $a_j\in \mathbb{G}_{a}\setminus \{0\}$, then we consider the product, 
\begin{align*}
x_{\beta_{k}}(-a_{1})x_{\beta_{2}}(-a_{2})\cdots x_{\beta_{1}}(-a_{k}) \xi,
\end{align*} 
which is denoted by $\xi'$. 
Clearly, $\xi'$ is in the $B_{L(w)}$-orbit of $\xi$, implying that $B_{L(w)}\xi'=B_{L(w)}\xi$.
Since $\dim({\rm Stab}_{B_{L(w)}}(\xi))=0$, we have $\dim({\rm Stab}_{B_{L(w)}}(\xi'))=0$.  
So, we may replace $\xi$ with the point $\xi'$.
\medskip

Since $P(w)$ is the stabilizer of $X_{wB}$ in $G$ and since the element $w_{0, J(w)}$ is represented by an element in the normalizer of $T$ in $L(w)\subset P(w)$, we see that $w_{0,J(w)}w$ represents a $T$-fixed point in $X_{wB}$. 
In particular, we see that $w_{0,J(w)}w \leq w$ in the Bruhat order of $W$. 
Equivalently, if $v$ denotes the element $w_{0,J(w)}w$, then we have the inclusion, 
\begin{align*}
X_{vB}\subset X_{wB}.
\end{align*}
We will compute the dimension of $X_{vB}$. 
Let $\gamma_{j}:=w_{0,J(w)}(\beta_{j})$ for $k+1\le j\le r$. 
Then note that $\gamma_{j}\in R^{+}$ for all $k+1\le j\le r$, and that $R^{+}(v^{-1})=\{\gamma_{j}: k+1\le j\le r\}=w_{0,J(w)}(R_{2})$. 
Thus, the dimension of $X_{vB}$ is given by $\# R_{2}=n=\dim T$. 
\medskip

Now, consider the point $\eta=w_{0,J(w)}\xi$ in $X_{wB}$. 
By replacing $w$ with $w_{0,J(w)}w_{0,J(w)} w$ (or, equivalently, replacing it with $w_{0,J(w)}v$) in (\ref{A:wlog}), we see that $\eta=\prod_{j=k+1}^{r}x_{\gamma_{j}}(b_{j})vB/B$ for some $b_{j}\in \mathbb{G}_{m}$. 
It follows that the point $\eta$ is in $UvB/B$. 
Consider the orbit map 
\begin{align*}
o_{\eta}:T\longrightarrow T\eta \subseteq X_{vB}: t\mapsto t\eta.
\end{align*}
We notice the following logical equivalences:
\begin{align*}
t\in {\rm Stab}_{T}(\eta) &\iff w_{0,J(w)}t w_{0,J(w)}\xi = \xi \\
&\iff w_{0,J(w)}t w_{0,J(w)} \in {\rm Stab}_{T}(\xi) \\
&\iff t \in w_{0,J(w)}{\rm Stab}_{T}(\xi)w_{0,J(w)}.
\end{align*}
After (\ref{A:wlog}), we observed that $\dim({\rm Stab}_{B_{L(w)}}(\xi))=0$. 
It follows that $\dim({\rm Stab}_{T}(\xi))=0$, or, equivalently, $\dim({\rm Stab}_{T}(\eta))=0$. 
But since $\dim(X_{vB})=\dim(T)$, $T\eta\subseteq X_{vB}$ is open in $X_{vB}$. 
Therefore, $X_{vB}$ is a toric variety. 
Hence, by the main result of \cite{Karuppuchamy}, $v$ is a Coxeter element. 
Therefore, we have $w=w_{0,J(w)}v$ and $\ell(w)=\ell(w_{0,J(w)})+\ell(v)$ as we claimed. 
\medskip

$(\Leftarrow)$ 
Conversely, since $c$ is a Coxeter element, by~\cite{Karuppuchamy}, $X_{cB}$ is a toric variety for the action of $T$. 
This means that there exists a point $\eta\in UcB/B$ such that $\overline{T\eta}=X_{cB}$.
Let $R^{+}(c^{-1})=\{\gamma_{j}: 1\le j\le n\}$. 
Then $\eta=\prod_{j=1}^{n}x_{\gamma_{j}}(a_{j})cB/B$ for some $a_{j}\in \mathbb{G}_{a}\setminus \{0\}$. 
We now consider the point $\xi=w_{0,J(w)}\eta$. 
Let $\beta_{j}=w_{0, J(w)}(\gamma_{j})$ for $1\le j\le n$. 
Since $w=w_{0,J(w)}c$ is such that $\ell(w)=\ell(w_{0,J(w)})+\ell(c)$, $\beta_{j}\in R^{+}$ for all $1\le j\le  n$.
So, we have  $\xi=\prod_{j=1}^{n}x_{\beta_{j}}(b_{j})wB/B$ for some $b_{j}\in \mathbb{G}_{m}$. 
Since $\dim({\rm Stab}_{T}(\eta))=0$ and $\ell(w)=\ell(w_{0,J(w)})+\ell(c)$, we have $\dim(B_{L(w)}\xi)=\dim(X_{wB})$. 
Therefore, it follows that $\overline{B_{L(w)}\xi}=X_{wB}$.  
\qed\end{proof}

For $J\subseteq J(w)$, let $L_{J}$ denote the standard Levi factor of $P_{J}$. 
Let $B_{J}$ denote the Borel subgroup of $L_{J}$. 
We are now ready to prove the first stated result of our paper, namely, Theorem~\ref{T:major}.
We paraphrase it for convenience.

\begin{Theorem}\label{thm6.2}
The Schubert variety $X_{wB}$ is a spherical $L_{J}$-variety such that $\dim X_{wB} = \dim B_{J}$ if and only if $w=w_{0,J}c$, where $w_{0,J}$ is the longest element of $W_{J}$ and $c$ is a Coxeter element of $W$ such that $\ell(w)=\ell(w_{0,J})+\ell(c)$.
\end{Theorem}

The proof of this theorem is similar but not identical to the proof of Theorem~\ref{thm6.2}. 
We will omit explanations of the arguments that are used in the previous proof. 	
	
\begin{proof}
$(\Rightarrow)$
$X_{wB}$ is a spherical $L_{J}$-variety such that $\dim X_{wB} = \dim B_{J}$. Then there exists a point $\xi$ in $X_{wB}$ such that $\overline{B_{J}\xi}=X_{wB}$ and $\dim({\rm Stab}_{B_{J}}(\xi))=0$. 
Let $R^{+}(w^{-1})=\{\beta_{1},\ldots, \beta_{r}\}$ be such that  $S_{1}=\mathbb{Z}J\cap R^{+}=\{\beta_{1},\ldots,\beta_{k}\}$. Let $S_2$ 
be the subset defined by 
\[
S_{2}:=\{\beta\in R^{+}(w^{-1}): supp(\beta)\nsubseteq J \}=R^{+}(w^{-1})\setminus S_{1}=\{\beta_{k+1},\ldots, \beta_{r}\}.
\] 
Note that $$\dim(X_{wB})=\# S_{1}+ \# S_{2}=\dim (B_{J})=\dim(T) + \dim(U_{J}),$$ where $U_{J}$ denotes the unipotent radical of $B_{J}$.  
Since $\dim(U_{J})=\ell(w_{0,J})=\# S_{1}$, by the above equation we have $\# S_{2}=\dim(T)$. Thus, we have $k=r-n$, where $n={\rm rank}(G)$.
\medskip

We observe as in the proof of Theorem~\ref{thm6.2} that, for $U^{*}=\{u\in U: x_{\beta}(u)\neq 0\ \text{ for all }\beta \in R^{+}\}$, 
the intersection $B_{L(w)}\xi\ \cap\ U^*wB/B\neq \emptyset$ is nonempty.
As before, we have $\xi=\prod_{j=1}^{r}x_{\beta_{j}}(a_{j})wB/B$ for some $a_{1},\ldots, a_{k}\in \mathbb{G}_{a}$ and $a_{k+1},\ldots, a_{r}\in \mathbb{G}_{a}\setminus \{0\}$. 
In fact, without loss of generality, we may assume that  $\xi=\prod_{j=k+1}^{r}x_{\beta_{j}}(a_{j})wB/B$ for some $a_{j}\in \mathbb{G}_{a}\setminus \{0\}$. 
Since $P_{J(w)}={\rm Stab}_{G}(X_{wB})$, the $T$-fixed point $w_{0,J}w$ satisfies $w_{0,J}w \le w$ in $W$. 
Let $v:=w_{0,J}w$. 
Let $\gamma_{j}=w_{0,J}(\beta_{j})$ for all $k+1\le j\le r$. 
Then note that $\gamma_{j}\in R^{+}$ for all $k+1\le j\le r$, and   $R^{+}(v^{-1})=\{\gamma_{j}: k+1\le j\le r\}=w_{0,J}(S_{2})$. 
Thus, $X_{vB}$ is a Schubert subvariety of $X_{wB}$ of dimension $\# R_{2}=n$. 
\medskip

We define $\eta:=w_{0,J}\xi \in X_{wB}$.
Notice that $\eta=\prod_{j=k+1}^{r}x_{\gamma_{j}}(b_{j})vB/B$ for some $b_{j}\in \mathbb{G}_{a}\setminus \{0\}$. 
This means that $\eta \in UvB/B$. 
Consider the orbit map $$o_{\eta}:T\longrightarrow T\eta \subseteq X_{vB}: t\mapsto t\eta.$$ 
Then the stabilizer ${\rm Stab}_{T}(\eta)=\bigcap_{j=k+1}^{r}\ker(\gamma_{j})$. 
So, we have $w_{0, J}({\rm Stab}_{T}(\eta))\subseteq {\rm Stab}_{T}(\xi)$. 
By the previous discussion we have $\dim({\rm Stab}_{T}(\eta))=0$. 
Since $\dim(X_{vB})=\dim(T)$, $T\eta\subseteq X_{vB}$ is open in $X_{vB}$. 
Therefore, $X_{vB}$ is a toric variety. 
Hence, by the main result of~\cite{Karuppuchamy}, $v$ is a Coxeter element. 
We now conclude that $w$ can be written in the form $w=w_{0,J}v$, where $\ell(w)=\ell(w_{0,J})+\ell(v)$.
\medskip

$(\Leftarrow)$ Conversely, since $c$ is a Coxeter element in $W$, by \cite{Karuppuchamy}, it follows that $X_{cB}$ is a toric variety for the action of $T$. 
So, there exists a point $\eta\in UcB/B$ such that $\overline{T\eta}=X_{cB}$.
Let $R^{+}(c^{-1})=\{\gamma_{j}: 1\le j\le n\}$. Then $\eta=\prod_{j=1}^{n}x_{\gamma_{j}}(a_{j})cB/B$ for some $a_{j}\in \mathbb{G}_{m}$. 
Consider the point $\xi=w_{0,J}\eta$. Let $\beta_{j}=w_{0, J}(\gamma_{j})$ for $1\le j\le n$. 
Since $w=w_{0,J(w)}c$ is such that $\ell(w)=\ell(w_{0,J(w)})+\ell(c)$, $\beta_{j}\in R^{+}$ for all $1\le j\le  n$.
So, we have  $\xi=\prod_{j=1}^{n}x_{\beta_{j}}(b_{j})wB/B$ for some $b_{j}\in \mathbb{G}_{m}$. 
Since $\dim({\rm Stab}_{T}(\eta))=0$ and $\ell(w)=\ell(w_{0,J})+\ell(c)$, we have $\dim(B_{J}\xi)=\dim(X_{wB})$. 
Therefore, it follows that $\overline{B_{J}\xi}=X_{wB}$.  
\qed\end{proof}

We demonstrate our Theorem~\ref{thm6.2} by an example. 
	
	\begin{Example}
	Let $w$ denote the permutation $5 13624$ in the symmetric group $\mathbf{S}_6$.  
	A reduced word of $w$ is given by 
	\begin{align*}
	513624 =  s_2 s_4 s_5 s_3 s_4 s_2 s_1.
	\end{align*} 
	Denoting the set $\{s_2,s_4\}$ by $J$, we see that $J\subseteq J(w)$. 
	The standard Levi subgroup of $GL(6,\C)$ associated with $J$ is given by 
	\[
	L_J:=GL(1,\C)\times GL(2,\C)\times GL(2,\C)\times GL(1,\C).
	\]
	Let $c$ denote $s_5 s_3 s_4 s_2 s_1$, which is a Coxeter element. 
	The product $s_2s_4$ is the maximal element of the parabolic subgroup of $\mathbf{S}_6$ corresponding to $J$.
	Since $w=w_{0,J}c$ and $\ell(w) = \ell(w_{0,J})+\ell(c)$, according to our Theorem~\ref{thm6.2}, the Schubert variety $X_{wB}$ is a spherical $L_J$-variety in $GL(6,\C)/B$. 
	Of course, this is in conformity with the main result of~\cite{gao2021classification}.
	\end{Example}

	\medskip
	
	We now have a remark about our Theorem~\ref{thm6.2} showing that it recovers the conjecture of Gao, Hodges and Yong~\cite{GHY} in characteristic 0.

	\begin{Remark}\label{R:relaxation}
	The assumption $\dim B_{L(w)}=\dim X_{wB}$ can be relaxed by accordingly relaxing the condition on the length of the (Coxeter) element $c$. 
	More precisely, the proof of Theorem~\ref{thm6.2} shows that the unipotent radical of $B_{L(w)}$ always acts faithfully on the open cell of $X_{wB}$. 
	Thus, if the inequality $\dim B_{L(w)} > \dim X_{wB}$ holds, then we see that a torus quotient of $B_{L(w)}$ acts faithfully with an open orbit.  
	In this case, we observe as in the proof of Theorem~\ref{thm6.2} that the orbit closure of the semisimple part of the torus quotient of $B_{L(w)}$ is a toric Schubert subvariety $X_{vB}\subset X_{wB}$ such that $w=w_{0,J}v$, where $J$ is a subset of $J(w)$.
	\end{Remark}

	\medskip

	We proceed to present an analogous sphericalness result for the BSDH-varieties. 
	To this end, let $\underline{w}$ be a reduced word in $S$ of the form $\underline{w}=(s_{i_1},s_{i_2},\ldots, s_{i_{r}})$. 
	Then we have the corresponding sequence of simple roots, $(\alpha_{i_1},\ldots, \alpha_{i_r})$. 
	Using this list, we define 
	\[
	J(\underline{w}):=\{\alpha_{i_{j}}: s_{i_{j}}s_{i_{k}}=s_{i_{k}}s_{i_{j}}  \text{ for all $1\le k\le j$} \}.
	\] 
	In other words, $J(\underline{w})$ is the set of simple roots $\alpha_{i_j}$ from the list $(\alpha_{i_1},\ldots, \alpha_{i_r})$ such that all of the simple roots $\alpha_{i_k}$ with $1\leq k \leq j$ commute with $s_{i_j}$.
	Let $P_{J(\underline{w})}$ be the standard parabolic subgroup of $G$ corresponding to the subset $J(\underline{w})$ of $S$. 
	Then there is a natural action of $P_{J(\underline{w})}$ on $X_{\underline{w}}$ that is given by the left multiplication.

\begin{Theorem}\label{T:sphericalBSDH}
Let $\underline{w}$ be a reduced word in $S$.
Then $X_{\underline{w}}$ is a spherical $L(\underline{w})$-variety such that $\dim B_{L(\underline{w})}=\dim X_{\underline{w}}$ if and only if $w_{0,J(\underline{w})}w$ is a Coxeter element, where $w$ is the element of $W$ associated with $\underline{w}$ and $w_{0,J(\underline{w})}$ denotes the longest element of $W_{J(\underline{w})}$.
\end{Theorem}

\begin{proof}
Since $\underline{w}$ is a reduced word, the natural product map 
	\begin{align*}
	\mathbf{m} :\ X_{\underline{w}} &\longrightarrow X_{wB} \\
	[p_1,\dots, p_m] & \longmapsto p_1\cdots p_mB, 
	\end{align*}
is a surjective birational morphism. 
Furthermore, since ${\mathbf{m}}$ is $P_{J(\underline{w})}$-equivariant, the BSDH-variety $X_{\underline{w}}$ is a spherical $L(\underline{w})$-variety 
if and only if $X_{wB}$ is a spherical $L(\underline{w})$-variety. 
But notice that $J(\underline{w}) \subseteq J(w)$.
Hence, it follows from Theorem \ref{thm6.2} that $X_{\underline{w}}$ is a spherical $L(\underline{w})$-variety such that
$\dim X_{\underline{w}} = \dim B_{L(\underline{w})}$ if and only if $w_{0, J(\underline{w})}w$ is a Coxeter element.

\qed\end{proof}

The following generalization of Theorem~\ref{T:sphericalBSDH} follows from our remark Remark~\ref{R:relaxation}.

\begin{Theorem}
Let $\underline{w}$ be a reduced word.
Then $X_{\underline{w}}$ is a spherical $L(\underline{w})$-variety if and only if $w_{0,J(\underline{w})}w$ is a product of distinct simple reflections, where $w$ is the element of $W$ associated with $\underline{w}$ and $w_{0,J(\underline{w})}$ denotes the longest element of $W_{J(\underline{w})}$.
\end{Theorem}

\begin{Corollary}\label{C:sphericalBSDH}
Let $\underline{w}$ be a reduced word.
Then, $X_{\underline{w}}$ is a spherical $L(s_{i_{1}})$-variety if and only if $s_{i_{1}}w$ is a product of distinct simple reflections.
\end{Corollary}

\section{Additional Results}\label{S:Applications}

There is an interesting interaction between the modality of the natural $T$-action and the complexity of a Schubert variety.
In~\cite[Theorem 1.4]{CanDiaz}, by using the main result of the article~\cite{Gaetz}, Can and Diaz obtained the following result.

\begin{Proposition}
Let $X_{wB}$ be a singular Schubert variety in $GL(n,\C)/B$. 
If the complexity of the $T$-action on $X_{wB}$ is 1, then $X_{wB}$ is a spherical $L_{J(w)}$-variety.
\end{Proposition}

Here we provide another proof.

\begin{proof}

In~\cite[Theorem 1.3, part 4]{LeeMasudaPark}, Lee, Masuda, and Park show that if $X_{wB}$ is a singular $T$-complexity 1 
Schubert variety in $GL(n,\C)/B$, then there is a reduced decomposition $\text{red}(w)$ of $w$ and an index $j\in \{1,\dots, n-1\}$ such that the segment $s_{j+1}s_j s_{j+2} s_{j+1}$ appears in $\text{red}(w)$ and no other simple reflection in $\text{red}(w)$ appears more than once (except $s_{j+1}$ in the segment $s_{j+1}s_j s_{j+2} s_{j+1}$). 
Now let us look at the product $s_{j+1}w$. 
Since both $s_j$ and $s_{j+2}$ appear on the right hand side of the left most $s_{j+1}$ in $\text{red}(w)$, the multiplication of $w$ on the left by $s_{j+1}$ deletes the leftmost $s_{j+1}$ from $\text{red}(w)$. 
But this means that $s_{j+1}w$ is a product of distinct simple reflections. 
Now, by Corollary~\ref{C:sphericalBSDH}, we see that $X_{\underline{w}}$ is a spherical $L(s_{j+1})$-variety.
Since the multiplication map $\mathbf{m} : X_{\underline{w}} \to X_{wB}$ is an $L(s_{j+1})$-equivariant birational morphism, 
we see that $X_{wB}$ is a spherical $L(s_{j+1})$-variety. 
At the same time, we know that $L(s_{j+1})$ is a subgroup of the stabilizer ${\rm Stab}_G(X_{wB})$.
Therefore, the bigger Levi subgroup $L_{J(w)}$ acts spherically on $X_{wB}$ as well. 
This finishes the proof of our assertion. 
\qed\end{proof}

We proceed with the assumption that $G$ is a connected semisimple algebraic group over an algebraically closed field $k$ of arbitrary characteristic. As before, for $w\in W$, let $P_{J(w)}$ denote the stabilizer of $X_{wB}$.
Let $L(w)$ denote the standard Levi factor of $P_{J(w)}$ so that $T\subseteq L(w)$.
Then the intersection $B\cap L(w),$ denoted $B_{L(w)},$ is a Borel subgroup of $L(w)$.	
For a subset $J\subset L(w)$, we follow the similar convention that $P_J$ denotes the parabolic subgroup generated by $B$ and $J$, $L_J$ denotes the standard Levi factor of $P_J$, and $B_{L_J}$ denotes the Borel subgroup $B\cap L_J$. 

We now fix two subsets $I$ and $K$ of $S$ such that $I\subset K$ and $K\neq \emptyset$.
For $w\in W^I$, let $X_{wP_I}$ denote the Schubert variety in $G/P_I$, where $P_I$ is the parabolic subgroup determined by $I$ and $B$. 
Following~\cite[pg. 34, paragraph 2]{RichmondSlofstra},
we define the {\em (right) parabolic decomposition of $w$ with respect to $K$} as the unique decomposition $w= vu$,  where $v\in W^K$ and $u\in W_K\cap W^I$. 
Closely related to the notion of a parabolic decomposition is the stronger notion of a {\em Billey-Postnikov decomposition} 
(abbreviated to a {\em $BP$-decomposition}).  
A parabolic decomposition $w=vu$ ($v\in W^K,\ u\in W^I\cap W_K$) is called a {\em BP-decomposition with respect to $(I,K)$} if the Poincar\'e polynomial of $X_{wP_I}$ is the product of the Poincar\'e polynomials of $X_{uP_I}$ and $X_{vP_K}$. 	
Note that if $I=\emptyset$, then we have $W^I= W$. 
We are now ready to prove our Theorem~\ref{T:inverseissmooth}.
We recall its statement for convenience. 	
\medskip

Let $w\in W$. 
Let $J$ be a subset of $J(w)$.
We assume that $X_{wB}$ is a spherical $L_J$-variety such that $\dim X_{wB}=\dim B_{L_J}$.  
Then the following assertions are equivalent: 
\begin{enumerate}
\item $X_{w B}$ is a smooth Schubert variety.
\item $X_{w^{-1}B}$ is a smooth Schubert variety.
\item $X_{c^{-1}P_J}$ is a smooth toric variety.
\end{enumerate}

\medskip
	
\begin{proof}[Proof of Theorem~\ref{T:inverseissmooth}]
(2 $\Leftrightarrow$ 3)
Under our assumptions, by Theorem~\ref{thm6.2}, we know that there exists an element $c\in W$, which is a product of distinct simple reflections, such that $w= w_{0,J}c$ and $\ell(w) = \ell(w_{0,J})+\ell(c)$. 
Clearly, the product $c^{-1}w_{0,J(w)}$ is the right parabolic decomposition of $w^{-1}$ with respect to $J$.
In fact, it follows from~\cite[Proposition 4.2 (c)]{RichmondSlofstra} that $c^{-1}w_{0,J}$ is a $BP$-decomposition with respect to $(\emptyset,J)$. 

Let $\pi : G/B\to G/P_{J}$ denote the canonical projection map. 
This is a smooth and proper morphism. 
On one hand, since $w^{-1} = c^{-1} w_{0,J}$ is the right parabolic decomposition of $w^{-1}$ with respect to $J$, we have $X_{c^{-1}P_{J}}= \pi( X_{w^{-1}B})$.
On the other hand, we know that $X_{c^{-1}B}$ is a toric variety in $G/B$ and $\pi (X_{c^{-1}B}) = X_{c^{-1}P_{J}}$. 
But at the same time, $X_{w^{-1}B}$ is the full preimage of $X_{c^{-1}P_{J}}$ under the morphism $\pi :G/B\to G/P_{J}$. 
Hence, the generic fiber of this surjection is given by $X_{w_{0,J}B}$, which is a smooth Schubert variety. 
Now, the rest of the proof of the equivalence of 2. and 3. follows from~\cite[Theorem 3.3, part 1.]{RichmondSlofstra}, which states that 
$X_{w^{-1}B}$ is smooth if and only if both the fiber and the base of the projection $X_{w^{-1}B}\to X_{c^{-1}P_J}$ are smooth Schubert varieties.

(1 $\Leftrightarrow$ 2) This result is proven by Carrell in~\cite[Corollary 4]{Carrell2011}.
\qed\end{proof}

\begin{Example}
Let us consider the Schubert variety $X_{wB}$ in $SL(4,\C)/B$, where $w$ is the element $w=s_1s_3s_2s_1s_3 = 4231$. 
Note that $w^2 = id$. We denote $J(w)$ simply by $J$. 
Since the longest element of $W_J=\{s_1,s_3\}$ is $w_{0,J}=s_1s_3$ and $c:=s_1s_3s_2$ is a Coxeter element, we see that $X_{wB}$ is a spherical $L_J$-variety such that $\dim B_{L_J}=\dim X_{wB}$.
It follows from Lakshmibai-Sandhya criterion~\cite{LakshmibaiSandhya} that $X_{w^{-1}B}$ is not a smooth Schubert variety. 
Conforming with our Theorem~\ref{T:inverseissmooth}, the toric variety $X_{c^{-1}P_{\{s_1,s_3\}}}$ is not smooth.
Indeed, $X_{c^{-1}P_{\{s_1,s_3\}}}$ is the unique singular Schubert variety in the Grassmann variety of two dimensional subspaces of $\C^4$.
\end{Example}

Our next observation is a consequence of the recent work~\cite[Proposition 4.9, p. 15]{NPS}.
It is concerned with the automorphism groups of the spherical $G$-Schubert varieties and the rigidity question for $G$-BSDH varieties. 

\begin{Theorem}
We assume that the underlying field of definitions is $k=\mathbb{C}$. 
Let $G$ be of simply-laced type. 
In addition, we assume that $G$ is simply connected. 
If $G\times_B X_{wB}$ is a spherical $G$-Schubert variety, then the identity component of the automorphism group of $G\times_B X_{wB}$ is $G_{ad}$,
that is, the adjoint group of $G$.
In addition, the first cohomology of its tangent sheaf vanishes.
\end{Theorem}

\begin{proof}
By Proposition~\ref{P:fourequivalent}, the $G$-Schubert variety $G\times_{B} X_{wB}$ is a spherical $G$-variety if and only if $X_{wB}$ is a toric variety. In particular, $G\times_{B} X_{wB}$ is isomorphic to the $G$-BSDH variety, $G\times_{B} X_{\underline{w}}$.
Hence, it is smooth. 
We know from~\cite[Proposition 4.9, p. 15]{NPS} that 
\begin{itemize}
	\item [(i)] ${\rm Aut}^{0}(G\times_B X_{wB})=G_{ad}$,
	\item [(ii)] $H^j(G\times_{B} X_{wB}, T_{G\times_{B} X_{wB}})=0$ for $j\ge 1$.
\end{itemize}
Here, ${\rm Aut}^{0}(G\times_{B}X_{wB})$ stands for the connected component of the identity element of the automorphism group of $G\times_BX_{wB}$, and $T_{G\times_{B} X_{wB}}$ denotes the tangent sheaf of $G\times_{B} X_{wB}$.
Hence, by~\cite[Proposition 6.2.10,p.272]{Huy}, we conclude that the spherical $G$-Schubert varieties are locally rigid for simply-laced groups. 
This finishes the proof of our assertion.
\qed\end{proof}

\begin{Remark}
If $G$ is not simply-laced, then spherical $G$-Schubert varieties need not be {\em rigid}, that is to say, $H^1(G\times_{B}X_{wB}, T_{G\times_{B} X_{wB}}) \neq 0$. 
We will demonstrate this phenomenon in our next example.
\end{Remark}

\begin{Example}\label{E:nonrigid}
Let $G=SO(5,\mathbb{C})$. 
Throughout this example, we work with a Borel subgroup $B\subset G$ that is associated with the negative roots. 
Let $w=s_{2}s_{1}$ and $w_{1}=s_{2}$. 
Let $p: X_{\underline{w}}\to X_{\underline{w_1}}$ denote the natural $P_{\alpha_{1}}/B$ $(\simeq \mathbb{P}^{1})$ fibration.
Then, we have the following short exact sequence of tangent sheaves on $X_{\underline{w}}$:
\begin{equation*}
0\longrightarrow \mathcal{L}_{\alpha_{1}}\longrightarrow T_{\underline{w}}\longrightarrow p^{*}T_{\underline{w_1}}\longrightarrow 0,
\end{equation*}
where $T_{\underline{w}}$ (respectively, $T_{\underline{w_1}}$) denotes the tangent sheaf of $X_{\underline{w}}$ (respectively, of $X_{\underline{w_1}}$), and $\mathcal{L}_{\alpha_1}$ denotes the relative tangent bundle with respect to the map $p$. 
Notice that $H^1(X_{\underline{w}}, \mathcal{L}_{\alpha_{1}})$ is isomorphic to the one dimensional $B$-module corresponding to the character $\alpha_{1}+\alpha_{2}$ of $B$, and $H^0(X_{\underline{w_1}}, T_{\underline{w_1}})_{\mu} = 0$ for $\mu = \alpha_1+\alpha_2.$ 
Therefore, by using the long exact sequence associated to the above short exact sequence, we obtain the following exact sequence:
\begin{equation*}
0 \longrightarrow H^1(X_{\underline{w}}, \mathcal{L}_{\alpha_{1}})\longrightarrow H^1(X_{\underline{w}}, T_{\underline{w}})\longrightarrow H^1(X_{\underline{w_1}}, T_{\underline{w_1}})\longrightarrow 0.
\end{equation*}
Note that $H^1(X_{\underline{w_1}},T_{\underline{w_1}} )$ vanishes.
Hence, we have the isomorphism,
\[
H^1(X_{\underline{w}}, T_{\underline{w}})= H^1(X_{\underline{w}},\mathcal{L}_{\alpha_{1}}).
\]
Moreover, $H^{1}(G\times_{B}X_{\underline{w}}, T_{G\times_{B} X_{\underline{w}}})$ fits in the following short exact sequence: 
\begin{equation*}
	0 \to H^1(G/B, \mathcal{L}(H^0(X_{\underline{w}},\mathcal{L}_{\alpha_{1}})))\to H^1(G\times_{B}X_{\underline{w}}, T_{G\times_{B}X_{\underline{w}}})\to H^0(G/B, \mathcal{L}(H^1(X_{\underline{w_1}},\mathcal{L}_{\alpha_{1}})))\to 0.
\end{equation*}
Since $\alpha_{1}+\alpha_{2}=\varpi_{1}$, we have $H^0(G/B, \mathcal{L}_{\alpha_{1}+\alpha_{2}})=V(\varpi_1)$, 
where $V(\varpi_1)$ is the Weyl module of highest weight $\varpi_1$.
Here, $\varpi_1$ is the fundamental weight associated with the simple root $\alpha_1$. 
Hence, we have $H^1(G\times_{B}X_{\underline{w}}, T_{G\times_{B}X_{\underline{w}}})\neq 0.$
Therefore, the $G$-Schubert variety $G\times_{B} X_{wB}$ is a not rigid $G$-variety. 
But we know from Proposition~\ref{P:fourequivalent} that it is a spherical $G$-variety.
\end{Example}

\section{Acknowledgements}
We express our gratitude to the referee for their meticulous review of our manuscript and for providing insightful suggestions that have significantly enhanced the quality of our paper.
We thank Ed Richmond for useful conversations about parabolic decompositions.
We thank Alex Yong for many useful conversations about Schubert varieties. 
Finally, we thank Yibo Gao, Reuven Hodges, and Alex Yong for sharing with us the details of their work. 
The first author thanks the Department of Mathematics at IIT Bombay for hospitality, where this collaboration has started. 
The first author was partially supported by a grant from the Louisiana Board of Regents (Contract no. LEQSF(2021-22)-ENH-DE-26).

\bibliography{references}
\bibliographystyle{plain}

\end{document}